\documentclass{gtart_a}
\pdfoutput=1
\usepackage{graphicx}

\title{Surface subgroups and handlebody attachment} 
\author{Vivien R Easson} 
\givenname{Vivien R}
\surname{Easson}
\address{Centre for Mathematical Sciences\\
University of Cambridge\\\newline
Wilberforce Road\\Cambridge CB3 0WB, UK} 
\email{V.Easson@dpmms.cam.ac.uk}
\urladdr{http://www.dpmms.cam.ac.uk/~vre20}

\volumenumber{10}
\issuenumber{}
\publicationyear{2006}
\papernumber{15}
\lognumber{0429}
\startpage{557}
\endpage{591}

\doi{}
\MR{}
\Zbl{}

\arxivreference{}
\arxivpassword{}   


\received{4 March 2004}
\revised{19 July 2005}
\accepted{25 March 2006}
\published{4 May 2006}
\publishedonline{4 May 2006}
\proposed{Cameron Gordon}
\seconded{Joan Birman, Colin Rourke}
\corresponding{}
\editor{}
\version{}

\arxivreference{}
\arxivpassword{}

\let\xysavmatrix\xymatrix
\def\xymatrix{\disablesubscriptcorrection\xysavmatrix}
\AtBeginDocument{\let\bar\wbar\let\tilde\wtilde\let\widetilde\wwtilde}

\makeatletter
\def\cnewtheorem#1[#2]#3{\newtheorem{#1}{#3}[section]
\expandafter\let\csname c@#1\endcsname\c@Thm}

\newtheorem{Thm}{Theorem}[section]
\newtheorem{Thm*}{Theorem}
\cnewtheorem{Lem}[Thm]{Lemma}
\cnewtheorem{Cor}[Thm]{Corollary}
\cnewtheorem{Prop}[Thm]{Proposition}
\theoremstyle{definition}
\cnewtheorem{Def}[Thm]{Definition}
\cnewtheorem{Property}[Thm]{Property}
\makeatother
\makeautorefname{Thm}{Theorem}
\makeautorefname{Thm*}{Theorem}
\makeautorefname{Lem}{Lemma}
\makeautorefname{Cor}{Corollary}
\makeautorefname{Def}{Definition}
\makeautorefname{Prop}{Proposition}

\newcommand{\Nb}{\ensuremath{\mathbb{N}}}
\newcommand{\Rb}{\ensuremath{\mathbb{R}}}
\newcommand{\Zb}{\ensuremath{\mathbb{Z}}}
\newcommand{\Hb}{\ensuremath{\mathbb{H}}}

\newcommand{\Psii}{\ensuremath{\Psi(\pa \widetilde{M}, \widetilde{\Sigma}_i)}}

\newcommand{\Pc}{\ensuremath{\mathcal{P}}}
\newcommand{\Lc}{\ensuremath{\mathcal{L}}}
\newcommand{\Wc}{\ensuremath{\mathcal{W}}}
\newcommand{\Hc}{\ensuremath{\mathcal{H}}}
\newcommand{\Gc}{\ensuremath{\mathcal{G}}}
\renewcommand{\C}{\ensuremath{\mathcal{C}}}

\newcommand{\Gm}{\ensuremath{\Gamma}}
\newcommand{\ra}{\ensuremath{\rightarrow}}
\newcommand{\al}{\ensuremath{\alpha}}
\newcommand{\sg}{\ensuremath{\sigma}}
\newcommand{\ep}{\ensuremath{\epsilon}}

\newcommand{\de}{\ensuremath{\delta}}
\newcommand{\pa}{\ensuremath{\partial}}

\newcommand{\Hpants}{\ensuremath{H}\negthinspace--\mbox{pants }}
\newcommand{\Ucollar}{\ensuremath{U}\negthinspace--\mbox{collar }}
\newcommand{\Uthin}{\ensuremath{U}\negthinspace--\mbox{thin }}

\begin{document}

\begin{abstract}
The main theorem of this paper generalizes recent results in Dehn
surgery to the case of handlebody attachment. We consider attaching
handlebodies and solid tori to the boundary of an irreducible,
boundary-irreducible, atoroidal and acylindrical 3--manifold. We show
that for a large class of homeomorphisms attaching these handlebodies,
the fundamental group of the resulting manifold contains the
fundamental group of a closed surface of genus at least two.
\end{abstract}

\maketitle


\section{Introduction}
\label{sec:intro}

In this paper, we discuss a particular generalization of the technique
of Dehn surgery for hyperbolic 3--manifolds, which we call handlebody
attachment. We consider a simple 3--manifold $M$ (not $B^3$) with
non-empty boundary, where \textsl{simple} means compact, irreducible,
$\pa$--irreducible, atoroidal and acylindrical.

As a consequence of the proof of Thurston's geometrization of Haken
manifolds, we can give $M$ a complete hyperbolic structure such that
the non-cuspidal part of $M$ is compact and every boundary component
of $M$ which has genus at least two inherits a totally geodesic metric
from the hyperbolic structure on $M$. When we use the term
\textsl{totally geodesic boundary} throughout this paper, this is what
we mean. Thus the only non-compact ends of $M$ are cusps.

We consider attaching a collection of handlebodies and solid tori
$\Hc$ to $\pa M$. Denote the resulting manifold by $M \cup_\phi
\Hc$, where $\phi=\{\phi_{\ell}\}$ is a collection of homeomorphisms,
and each $\phi_{\ell}$ is a map from a boundary component $(\pa
M)_\ell$ of $\pa M$ to the boundary of a handlebody or solid torus
$H_{\ell}$ in $\Hc$.

A \textsl{surface subgroup} of $M$ is a subgroup $\Gm \leq \pi_1(M)$
which is isomorphic to the fundamental group of a closed surface of
genus at least two. In~\cite{CooperLongReid:1997}, Cooper--Long--Reid
showed that $M$ always has a non-peripheral surface subgroup.

In the Dehn surgery case, $\Hc$ is a collection of solid tori. Recent
papers of Cooper--Long~\cite{CooperLong:2001} and Li~\cite{Li:2002}
showed independently that all but finitely many Dehn surgeries on a
one-cusped hyperbolic manifold give rise to a manifold with a surface
subgroup. Bart~\cite{Bart:2001} gives a result when there are many
cusps. We generalize these results to the case of handlebody
attachments.

For a handlebody $H \in \Hc$, we consider attachments involving
pseudo-Anosov homeomorphisms $h \co \pa H \ra \pa H$ whose stable
lamination is of full type with respect to some pants decomposition of
$H$ along meridians, see \fullref{def:fulltype}. This condition
is generically true, and is equivalent (see~\cite{AbramsSch:2003}) to
the condition that $h$ has stable lamination which lies in the Masur
domain of the handlebody.

\begin{Thm*}[Theorems~\ref{thm:sfsubgp} and~\ref{thm:sfsubgp_many}]
\label{thm:intro_sfsubgp}
Suppose $M \neq B^3$ is a simple 3--manifold with $m+m' \geq 1$
boundary components, precisely $m'$ of which are tori. Take a
collection $\Hc=\{H_1,\ldots,H_{m+m'}\}$ of handlebodies and solid
tori whose genera match those of $\pa M$. Let $M \cup_\phi \Hc$ denote
the closed 3--manifold obtained by gluing each boundary component
$(\pa M)_{\ell}$ to $\pa H_{\ell}$ by a homeomorphism $\phi_\ell$.

Suppose moreover that $h_\ell \co \pa H_\ell \ra \pa H_\ell$ is a
homeomorphism which is either a pseudo-Anosov homeomorphism whose
stable lamination is of full type, or an Anosov homeomorphism,
according to whether $H_\ell$ has genus at least two or is a solid
torus respectively.

Given homeomorphisms $\phi'_\ell \co (\pa M)_\ell \to \pa H_\ell$,
there exist integers $(N_\ell)_{\textrm{min}}$ such that if $\phi_\ell
= h_\ell^{N_\ell} \circ \phi'_\ell$ with $N_\ell \geq
(N_\ell)_{\textrm{min}}$ for all $\ell$, the group $\pi_1(M
\cup_{\phi} \Hc)$ contains a surface subgroup.
\end{Thm*}

There are three main ingredients required in our proof of
\fullref{thm:intro_sfsubgp}. Firstly, a closed $\pi_1$--injective
surface $\Sigma$ with a finite cover which lifts to an embedded
incompressible surface in a finite cover of $M$. Secondly, suitable
classes of gluing homeomorphisms. Finally, a geometrical analysis of
the ways in which $\Sigma$ might fail to remain $\pi_1$--injective
under handlebody attachment.

We first prove \fullref{thm:intro_sfsubgp} assuming $M$ has no
toral boundary components: every component of $\pa M$ is a surface of
genus at least two. This proof is given in
Sections~\ref{sec:projwavelike} to~\ref{sec:incomp}, and generalized
to the case where $M$ also has some toral boundary components in
\fullref{sec:many_components}. The precise method of proof when
the chosen hyperbolic metric on $M$ has cusps depends on whether
$\Sigma$ has accidental parabolics, but is broadly similar to the
non-cusped case.

In both cases, we prove \fullref{thm:intro_sfsubgp} as a corollary
of the following theorem. 

\begin{Thm*}[Theorems~\ref{thm:sfs_incomp} and~\ref{thm:sfs_incomp_many}]
\label{thm:intro_sfs_incomp}
Let $M$ be a hyperbolic 3--manifold with $m$ totally geodesic boundary
components and $m'$ cusps, $m+m' \geq 1$. Let $\Sigma$ be a connected,
orientable, closed, immersed surface in $M$ with a finite cover which
lifts to a non-peripheral embedded incompressible surface in a finite
regular cover $\pi \co \widetilde{M} \to M$ of degree $d$.

As in \fullref{thm:intro_sfsubgp}, let $M \cup_\phi \Hc$ denote
the closed 3--manifold obtained by attaching handlebodies to $\pa M$
via homeomorphisms $\phi_{\ell} \co (\pa M)_\ell \ra \pa H_\ell$ and
consider maps $h_\ell$ satisfying the condition given in that
theorem. 

Given homeomorphisms $\phi'_\ell \co (\pa M)_\ell \to \pa H_\ell$,
there exist integers $(N_\ell)_{\textrm{min}}$ such that if $\phi_\ell
= h_\ell^{N_\ell} \circ \phi'_\ell$ with $N_\ell \geq
(N_\ell)_{\textrm{min}}$ for all $\ell$, the surface $\Sigma$ remains
$\pi_1$--injective in the resulting manifold $M \cup_\phi \Hc$.
\end{Thm*}

The existence of such a surface $\Sigma$ is given by
Cooper--Long--Reid~\cite{CooperLongReid:1997}. However our proof works
for any virtually embedded $\pi_1$--injective surface satisfying the
conditions above, not just those provided by their construction.

The second ingredient is making a good choice of attaching maps. The
crucial property possessed by each gluing homeomorphism described
above is that its dynamics are governed by that of the pseudo-Anosov
map. In \fullref{sec:projwavelike} we explain the condition for a
lamination to be of full type, in terms of the way it intersects pants
decompositions of the handlebody. 

When we attach handlebodies by maps as above, we force each member of
a certain finite collection of curves to have high geometric
intersection number with every curve bounding a disc in a handlebody
$H$ under the attachment map. Furthermore, we can use the geometry of
these intersection sets to prove a key lemma,
\fullref{prop:masurext}.

Now we use our final ingredient. Suppose that our surface $\Sigma$
fails to inject at the level of the fundamental group into $M
\cup_\phi \Hc$. Then there exists an essential loop $\Lc$ in $\Sigma$,
and an immersed disc $D$ in $M \cup_\phi \Hc$ which spans $\Lc$. The
intersection of $D$ with $M$ is a planar surface, and we may consider
an immersed least-area representative $Q$ of this planar surface. By
the Gauss--Bonnet theorem, the area of $Q$ is bounded above by a
linear function of its Euler characteristic, and all but one component
of $\pa Q$ bound a disc in $\Hc$ under the attaching maps.

Our proof relies on showing that the intersections as above contribute
area to a surface $R$. This surface is obtained as a cover of $Q$ of
degree at most $d$, where $d$ is the degree of a finite cover $\pi \co
\widetilde{M} \to M$ in which $\Sigma$ becomes embedded. For a
sufficiently high power of the pseudo-Anosov map, the intersections
contribute too much area, and we obtain a contradiction with the
Gauss--Bonnet theorem.

We develop this intuition in Sections~\ref{sec:projwavelike}
to~\ref{sec:incomp} in order to prove
\fullref{thm:intro_sfs_incomp} and hence
\fullref{thm:intro_sfsubgp} in the case when every boundary
component of $M$ is of genus at least two.  In
\fullref{sec:many_components} we show that the presence of cusps
does not materially affect our argument, so that
Theorems~\ref{thm:intro_sfsubgp} and~\ref{thm:intro_sfs_incomp} remain
true in the general case.

Note that attaching a solid torus by a map as in the theorems ensures
that the corresponding Dehn filling slope lies outside the union of a
finite number of infinite strips. In fact a similar argument to ours
shows that our results hold in the latter case as well. The important
condition turns out to be that the Dehn filling slope should have high
intersection number with a finite number of exceptional slopes. These
slopes depend only on the surface $\Sigma$ and on the way in which its
convex hull intersects the boundary of $M$.

I would like to thank my DPhil supervisor Marc Lackenby for a great
many helpful conversations. My thanks also to Juan Souto for helping
me to place my original condition on the pseudo-Anosov maps in the
context of the literature, to Saul Schleimer for comments on
terminology, and to the referee for several helpful suggestions.

\section{Defining the projective wavelike set}

\label{sec:projwavelike}

Let $H$ be a handlebody of genus $g$, and write $B(\pa
H)_{\textrm{emb}}$ for the collection of non-trivial embedded
1--manifolds on $\pa H$, each component of which bounds a disc in $H$,
considered up to isotopy. Similarly, let $B(\pa H)$ denote the
collection of non-trivial regular homotopy classes of closed curves on
$\pa H$, not necessarily embedded, each component of which bounds a
disc in $H$.

We will consider the first set as a subset of $ML(\pa H)$, measured
laminations on $\pa H$, by giving each curve component unit
weight. The second set may be thought of in terms of Bonahon's theory
of geodesic currents (e.g.~\cite{Bonahon:1988}), but this is not
necessary for our discussion.

\begin{Def}
\label{def:hpd}
Suppose $\Pc$ is a pants decomposition of the surface $\pa H$ such
that every boundary component of a pair of pants $P \in \Pc$ lies in
$B(\pa H)_{\textrm{emb}}$. Then we say that $\Pc$ is a
\textsl{handlebody pants decomposition} and that every $P \in \Pc$ is
an \textsl{$H$\negthinspace--pants}.
\end{Def}

We will find it useful to distinguish two different kinds of essential
subarcs in an \Hpants $P$, and relate this to the way that a
lamination $\lambda$ intersects $P$.

\begin{Def}
\label{def:waves}
Suppose that $P$ is a pair of pants, and that $\xi$ is an essential
subarc properly embedded in $P$. If both of the endpoints of $\xi$ lie
on the same boundary component of $P$, we say that $\xi$ is a
\textsl{wave} in $P$. Otherwise, $\xi$ connects two different boundary
components of $P$, and we say that $\xi$ is a \textsl{seam} in $P$.

Similarly, an essential immersed subarc $\xi$ with $\xi \cap \pa P =
\pa \xi$ is said to be an immersed wave or an immersed seam according
to whether its endpoints lie on the same or different components of
$\pa P$.
\end{Def}

\begin{figure}[ht!]\small
\begin{center}
\begin{picture}(0,0)%
\includegraphics{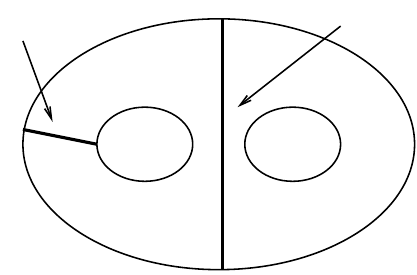}%
\end{picture}%
\setlength{\unitlength}{3947sp}%
\begin{picture}(1998,1304)(1276,-1590)
\put(2851,-361){\makebox(0,0)[lb]{\smash{wave}}}
\put(1951,-1411){\makebox(0,0)[lb]{\smash{$P$}}}
\put(1276,-436){\makebox(0,0)[lb]{\smash{seam}}}
\end{picture}%
\caption{Waves and seams in a pair of pants}
\end{center}
\end{figure}
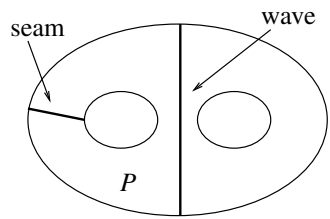

\begin{Def}
\label{def:triangleineq}
Consider a measured lamination $\lambda$ and an \Hpants $P$. Label the
boundary components of $P$ by $\pa P_1, \pa P_2, \pa P_3$.  Let
$a,b,c$ denote respectively $a=|\lambda \thinspace \cap \thinspace \pa
P_1|$, $b=|\lambda \thinspace \cap \thinspace \pa P_2|$ and
$c=|\lambda \thinspace \cap \thinspace \pa P_3|$, where $|x \thinspace
\cap \thinspace y|$ denotes the geometric intersection number of
measured laminations $x$ and $y$.

Then $\lambda$ satisfies \textsl{all triangle inequalities} on $P$ if
$a+b \geq c$, $b+c \geq a$ and $c+a \geq b$. We say that $\lambda$
satisfies \textsl{all strict triangle inequalities} if every such
inequality is strict.

On the other hand, if at least one of these strict inequalities does
not hold (e.g.~$a+b \leq c$), we say that $\lambda$ satisfies
\textsl{some triangle subequality}. Thus it is possible (take $a+b=c$,
$b+c \geq a$, $c+a \geq b$) that $\lambda$ satisfies both all triangle
inequalities and some triangle subequality.
\end{Def}

We can show that curves in $B(\pa H)$ usually have an immersed wave
with respect to some component of a handlebody pants decomposition. 

\begin{Lem}
\label{lem:bdhwave}
Consider a closed, connected, possibly non-embedded curve $\al$ which
lies in the boundary of a handlebody $H$ and is homotopically trivial
in $H$. Let $\Pc$ be a handlebody pants decomposition of $\pa H$. Then
either $\al$ is homotopic into some component $P \in \Pc$ or $\al$ has
an immersed wave with respect to some pants $P \in \Pc$.
\end{Lem}

\proof Suppose some component of $\al$ is not homotopic into
any component $P \in \Pc$, and write $f \co D \looparrowright H$ for
an immersed disc in $H$ that it bounds. Then by ambient isotopy of
$f(D)$ we may assume it intersects the discs in $\Pc$ in a minimal
non-empty collection of arcs whose preimages are properly and
disjointly embedded in $D$. An outermost such preimage arc $c$
separates off a subdisc $f \co D' \looparrowright H$ with boundary $c
\cup \al'$ whose image lies in the union of some \Hpants $P$ and the
discs in $H$ bounded by its boundary components.

\begin{figure}[ht!]\small
\begin{center}
\begin{picture}(0,0)%
\includegraphics{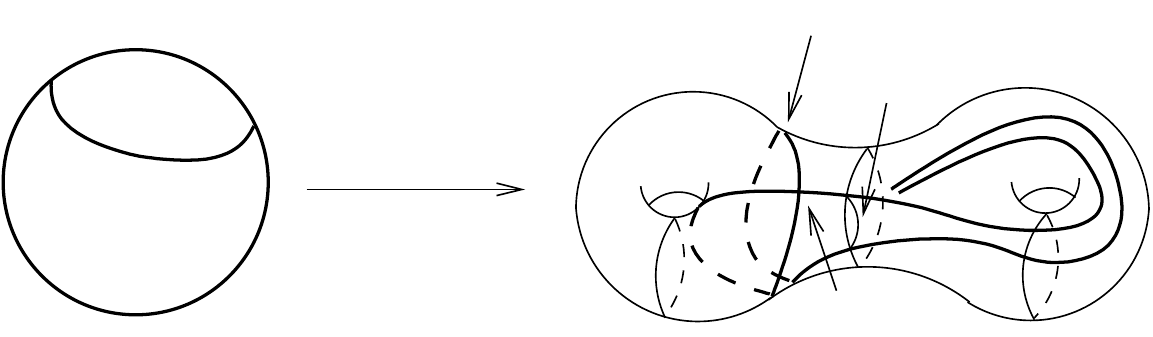}%
\end{picture}%
\setlength{\unitlength}{4144sp}%
\begin{picture}(5262,1539)(237,-995)
\put(715,-661){\makebox(0,0)[lb]{\smash{$D$}}}
\put(786,-18){\makebox(0,0)[lb]{\smash{$D'$}}}
\put(526,-203){\makebox(0,0)[lb]{\smash{$c$}}}
\put(1283,182){\makebox(0,0)[lb]{\smash{$\alpha'$}}}
\put(2001,-161){\makebox(0,0)[lb]{\smash{$f$}}}
\put(3216,-139){\makebox(0,0)[lb]{\smash{$P$}}}
\put(3908,-950){\makebox(0,0)[lb]{\smash{$f(D')$}}}
\put(4211, 96){\makebox(0,0)[lb]{\smash{$f(c)$}}}
\put(3837,424){\makebox(0,0)[lb]{\smash{$f(\alpha')$}}}
\end{picture}
\caption{A wave in a curve bounding a disc} 
\end{center}
\end{figure}
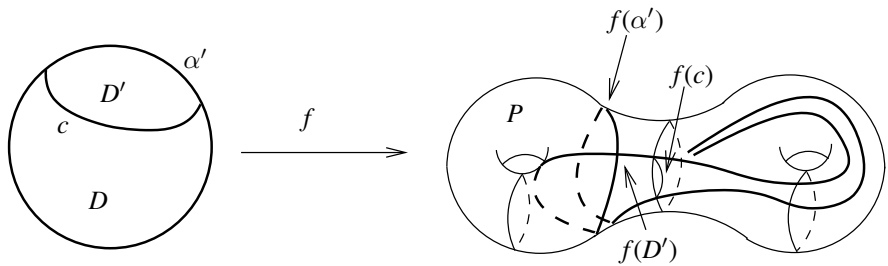

Since the endpoints of $f(c)$ lie in the same boundary component of
$P$, so do those of $f(\al')$. By minimality, $f(\al')$ is an
essential arc which is proper and immersed in $P$, so it forms an
immersed wave in $P$.\endproof

For $\al$ embedded, this wave will be embedded, and we can say more.

\begin{Cor}
Suppose that $\al \in B(\pa H)_{\textrm{emb}}$. Then the following
holds for any handlebody pants decomposition $\Pc$, whether or not
there is a wave. There exist an \Hpants $P \in \Pc$, and two boundary
components $\pa P_1$ and $\pa P_2$ of $P$, such that $\al$ has no
subarcs in $P$ which are seams connecting $\pa P_1$ and $\pa P_2$.
\end{Cor}

\proof If $\al$ has a wave, it cannot have some seam which
would cross that wave; and if $\al$ lies inside $P$, it has a
component parallel to a boundary component of some $P$ and so cannot
have any seam running from that component.\endproof

Moreover, suppose $\al \in B(\pa H)_{\textrm{emb}}$ has $k>0$ waves in
$P$ joining $\pa P_3$ to itself. These also separate $\pa P_1$ from
$\pa P_2$, so we deduce that the corresponding triple of intersection
numbers is of the form $(a,b,c)=(a,b,a+b+2k)$. If $\al$ has no waves
in $P$ but some component of $\al$ is parallel to $\pa P_1$, we have
$(a,b,c)=(0,b,b)$.

In particular, given any handlebody pants decomposition $\Pc$, and a
curve $\al$ in the set $B(\pa H)_{\textrm{emb}}$, there exists an
\Hpants $P \in \Pc$ such that $\al$ satisfies a triangle subequality
on $P$. We are interested in all curves and laminations satisfying
this conclusion, which motivates the following definition.

\begin{Def}
\label{def:wavelike}
Define the \textsl{wavelike set} $\Wc(H) \subset ML(\pa H)$ to be
\[ \Wc(H) = \{ \lambda : \forall \Pc, \exists P \in
\Pc \mbox{ s.t. } \lambda \mbox{ satisfies some triangle subequality on
$P$ }\}\] where $\Pc$ ranges over all handlebody pants decompositions
of $H$.
\end{Def}

Since a measured lamination $\lambda$ lies in $\Wc(H)$ if and only if
any positive scalar multiple of $\lambda$ does, we can also define the
\textsl{projective wavelike set} $P\Wc(H)$ consisting of all
projective measured laminations which have a representative in
$\Wc(H)$. This is a subset of $PML(\pa H)$, projective measured
lamination space. We may therefore make the following definition.

\begin{Def}
\label{def:fulltype}
Let $\lambda$ be a measured lamination in $ML(\pa H)$. If it does not
lie in $\Wc(H)$ then we say it is \textsl{of full type}. Similarly, a
projective measured lamination is of full type if it does not lie in
$P\Wc(H)$.

Suppose $h \co \pa H \to \pa H$ is a pseudo-Anosov homeomorphism whose
stable and unstable laminations $\lambda^S$ do not lie in the
projective wavelike set $P\Wc(H)$. Then we say that $h$ is \textsl{of
full type}.
\end{Def}

We showed above that $B(\pa H)_{\textrm{emb}} \subset
\Wc(H)$. Abrams--Schleimer proved that $P\Wc(H)$ is equal to the
closure of $PB(\pa H)_{\textrm{emb}}$ in the subspace of $PML(\pa H)$
consisting of those laminations which occur as the fixed points of
pseudo-Anosov homeomorphisms (see~\cite{AbramsSch:2003}, 11.5). They
also demonstrated that $h$ has stable lamination of full type if and
only if $h$ has stable lamination contained in the Masur domain of
$H$: the set of laminations having non-zero geometric intersection
with every lamination in the closure of $PB(\pa H)_{\textrm{emb}}$.

By work of Kerckhoff~\cite{Kerckhoff:1990} which extends results of
Masur~\cite{Masur:1986}, the closure of $PB(\pa H)_{\textrm{emb}}$ and
hence $P\Wc(H)$ has measure zero in $PML(\pa H)$. Thus a generic
pseudo-Anosov map is of full type. The terminology \emph{full type}
comes from Kobayashi \cite{Koba:1988pa}.

We will consider attaching maps of the form $\phi=h^N \circ \phi'$
which are formed as the composition of a power of a pseudo-Anosov
homeomorphism $h \co \pa H \to \pa H$ whose stable lamination is of
full type and a homeomorphism $\phi' \co(\pa M)_\ell \to \pa H$. Our
results hold for all $\phi$ with $N$ larger than some constant
depending on $h$ and $\phi'$.

We use throughout the conventions for stable and unstable laminations
of a pseudo-Anosov map found in~\cite{CassonBleiler:1988}.  In
\fullref{sec:crossing}, we show how to put these definitions
together with that of the set of distinguished curves defined below.

\section{The set of distinguished curves $X_i$}

\label{sec:cutlocus}

We start by defining what it means for a 3--manifold to be
simple. Then, given such a manifold $M$ (with non-empty boundary) and
a suitable $\pi_1$--injective surface $\Sigma$ immersed in it, we
define sets $X_i$ of distinguished curves which will play a crucial
role in our argument.

\begin{Def}
An orientable 3--manifold $M$ is said to be \textsl{simple} if it is
compact, irreducible, $\pa$--irreducible, atoroidal and acylindrical.
\end{Def}

Take a simple 3--manifold $M$ with $m$ boundary components of genus at
least two and $m'$ toral boundary components, $m+m' \geq 1$. As a
consequence of the proof of Thurston's geometrization theorem for
Haken 3--manifolds, we can give $M$ (minus its toral boundary) a
complete hyperbolic structure $g$ such that the non-cuspidal part of
$M$ is compact and every boundary component of genus at least two
inherits a totally geodesic hyperbolic metric from it.

For $M$ as above, we will consider an immersed surface $\eta \co
\Sigma \looparrowright M$ which is is a closed, connected, orientable
$\pi_1$--injective surface of genus at least two and which has a
finite cover $\widetilde{\Sigma}$ lifting to a non-peripheral
embedded incompressible surface in a degree $d$ cover $\pi \co
\widetilde{M} \to M$.

We can always assume that this cover is regular, by lifting if
necessary to a further finite cover. The next result follows from
various facts originally due to Thurston, as described below.

\begin{Lem}
Let $\eta \co \Sigma \looparrowright M$ be a surface in a simple
3--manifold $M$ as above. Then $\Sigma$ is geometrically finite.
\end{Lem}

\proof Suppose firstly that $M$ has at least one
boundary component of genus at least two. Then, giving $M$ a
hyperbolic metric $g$ as above, the limit set of the Kleinian group
$\pi_1(M)$ in the hyperbolic sphere at infinity $S^2_\infty = \pa
\overline{\Hb^3}$ does not fill up all of the
sphere. Matsuzaki--Taniguchi~\cite{MatsuTani:1998} call this a
Kleinian group of the second kind. Theorem 3.11
of~\cite{MatsuTani:1998} states that every finitely generated subgroup
of such a group is geometrically finite. This fact of Thurston's can
also be found in~\cite{MorganBass:1984}. Since $\Sigma$ is a closed
surface, $\pi_1(\Sigma)$ is finitely generated and hence geometrically
finite.

On the other hand, suppose that $M$ only has toral boundary
components, so that $(M,g)$ is a complete hyperbolic 3--manifold with
finite volume. In this case, it follows from work of
Thurston~\cite{Thurston:1979} and Bonahon~\cite{Bonahon:1986} that
either $\Sigma$ is geometrically finite or it is a virtual fibre of
$M$: see for example the discussion in
Menasco--Reid~\cite{MenascoReid:1992}. But $\Sigma$ is closed so $M$
cannot virtually fibre over it. Therefore $\Sigma$ is geometrically
finite.\endproof

Take a virtually embedded $\pi_1$--injective surface $\eta \co \Sigma
\looparrowright M$ as above. Write $\widetilde{\Sigma}_1, \ldots,
\widetilde{\Sigma}_s$ for the images of the connected surface
$\widetilde{\Sigma}$ which is a cover of $\Sigma$ lifting to
$\widetilde{M}$. These are embedded incompressible surfaces which
intersect each other in $\widetilde{M}$. We shall assume these
conditions and this terminology throughout this paper. Until
\fullref{sec:many_components}, we also assume $m'=0$, so $M$ has
no cusps.

Let $\Psi_i$ denote the characteristic submanifold of
$\widetilde{M}^\ast_i = \widetilde{M}-
\mbox{int}(N(\widetilde{\Sigma}_i))$, where $\widetilde{M}^\ast_i$ is
the manifold obtained from $\widetilde{M}$ by removing a small open
neighbourhood of the embedded surface $\widetilde{\Sigma}_i$. Its
vertical boundary $\pa_v \Psi_i$ is a collection of disjoint essential
annuli properly embedded in $(\widetilde{M}^\ast_i,\pa \widetilde{M}
\cup \widetilde{\Sigma}_i)$.

Moreover, since $M$ is hyperbolic, $\Psi_i$ consists of $I$--bundles
over subsurfaces of the boundary and Seifert-fibred solid tori. This
follows since $M$ being atoroidal implies that no other Seifert-fibred
spaces occur as components of $\Psi_i$. Again, we note that the
submanifolds $\Psi_1, \ldots, \Psi_s$ intersect each other.

Consider a restricted subset of the connected components of each
$\Psi_i$. This subset consists of those $I$--bundles whose vertical
boundary is a collection of annuli with one boundary component in each
of $\pa \widetilde{M}$ and $\widetilde{\Sigma}_i$. Write $\Psii$ for
the collection of all such $I$--bundles in $\widetilde{M}$.

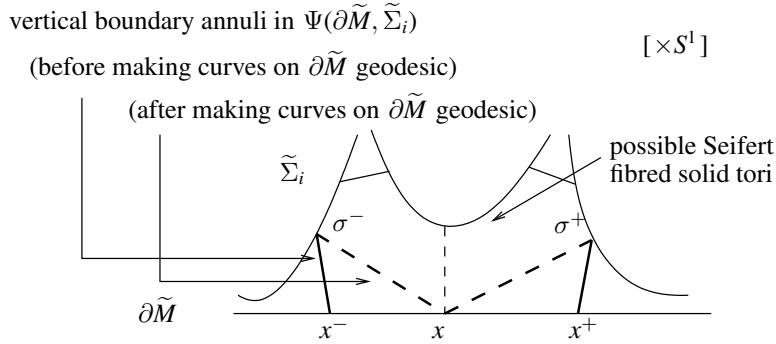
\begin{figure}[ht!]\small
\begin{center}
\begin{picture}(0,0)%
\includegraphics{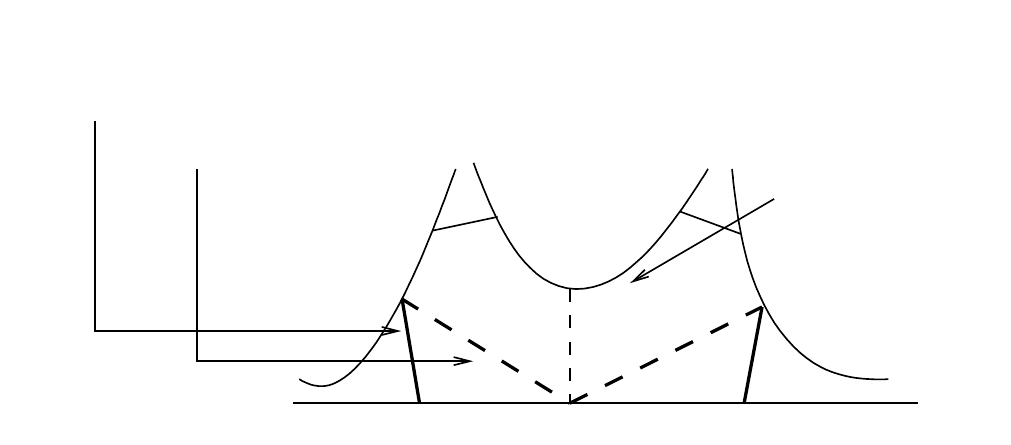}%
\end{picture}%
\setlength{\unitlength}{4144sp}%
\begin{picture}(4608,2025)(286,-2216)
\put(1891,-1231){\makebox(0,0)[lb]{\smash{$\widetilde{\Sigma}_i$}}}
\put(3624,-2176){\makebox(0,0)[lb]{\smash{$x^+$}}}
\put(2805,-2176){\makebox(0,0)[lb]{\smash{$x$}}}
\put(1044,-2078){\makebox(0,0)[lb]{\smash{$\partial \widetilde{M}$}}}
\put(406,-601){\makebox(0,0)[lb]{\smash{(before making curves on $\partial \widetilde{M}$ geodesic)}}}
\put(998,-856){\makebox(0,0)[lb]{\smash{(after making curves on $\partial \widetilde{M}$ geodesic)}}}
\put(286,-323){\makebox(0,0)[lb]{\smash{vertical boundary annuli in $\Psi(\partial \widetilde{M},\widetilde{\Sigma}_i)$}}}
\put(2206,-1531){\makebox(0,0)[lb]{\smash{$\sigma^-$}}}
\put(3879,-1074){\makebox(0,0)[lb]{\smash{possible Seifert}}}
\put(3879,-1239){\makebox(0,0)[lb]{\smash{fibred solid tori}}}
\put(4073,-466){\makebox(0,0)[lb]{\smash{[$\times S^1$]}}}
\put(3533,-1546){\makebox(0,0)[lb]{\smash{$\sigma^+$}}}
\put(2131,-2176){\makebox(0,0)[lb]{\smash{$x^-$}}}
\end{picture}%
\caption{The curve $x$ is a cut curve in this picture $\times$ $S^1$}
\label{fig:illustrating_cut_curves}
\end{center}
\end{figure}

The characteristic submanifold is a purely topological object, but we
will impose some geometrical restrictions on the subset $\Psii$. We
will assume that the boundary components of $\Psii \cap \pa
\widetilde{M}$ are geodesic curves in $\pa \widetilde{M}$. This may
make some previously disjoint annuli in $\pa_v \Psii$ intersect along
a closed geodesic in $\pa \widetilde{M}$, as shown in
\fullref{fig:illustrating_cut_curves}. Note that there are finitely
many of these curves since each of them is the geodesic straightening
of a boundary curve of the original $\Psii \cap \pa \widetilde{M}$. We
now give a precise definition.

\begin{Def}
\label{def:cutcurves}
Suppose $x$ is a closed simple geodesic on a totally geodesic boundary
component of $\pa \widetilde{M}$ with the following property. There
exist curves $\sg^+$ and $\sg^-$ on $\widetilde{\Sigma}_i$ which are
freely homotopic to $x$ in $\widetilde{M}^\ast_i$ but not freely
homotopic to each other in $\widetilde{\Sigma}_i$. Then we say that
$x$ is a \textsl{cut curve} for $\Psii$.
\end{Def}

We are interested in the intersection of $\Psii$ with the boundary
$\pa \widetilde{M}$. For our argument it is convenient to express this
intersection set in terms of the convex hull
$\C(\widetilde{\Sigma}_i)$ of $\widetilde{\Sigma}_i$. We will use the
same notation $\C(-)$ throughout for convex hulls in $\Hb^3$ and their
images when we project to a quotient manifold.

We will also occasionally find it useful to consider a slightly larger
submanifold containing the convex hull which has smooth boundary. We
may do this by extending the convex hull in a neighbourhood of small
width which will be determined precisely later. To avoid excessive
notation, we also let $\C(-)$ denote the resulting extended convex
hull with smooth boundary.

\begin{Prop}
\label{prop:hull_equality}
For $\widetilde{M}$, $\C(\widetilde{\Sigma}_i)$ and $\Psii$ as defined
above,
\[\Psii \cap \pa \widetilde{M} 
  = \C(\widetilde{\Sigma}_i) \cap \pa \widetilde{M}.\]
\end{Prop}

\proof For each boundary component $(\pa \widetilde{M})_\ell$
of $\widetilde{M}$, choose a basepoint $\ast$ on $(\pa
\widetilde{M})_\ell$ for all fundamental groups considered. Starting
from the right-hand-side of the above equation, we have
\begin{eqnarray*}
 \C(\widetilde{\Sigma}_i) \cap \pa \widetilde{M} &=& 
    \bigcup_\ell \left ( \C(\widetilde{\Sigma}_i) 
                     \cap (\pa \widetilde{M})_\ell 
              \right ) \\  
&=& \bigcup_\ell \left ( \C \left ( \Lambda(\pi_1(\widetilde{\Sigma}_i)) 
                         \right ) \cap \mbox{ }
                      \C \left 
                      (\Lambda(\pi_1((\pa \widetilde{M})_\ell)) 
                      \right )
              \right ),
\end{eqnarray*}
where $\Lambda(\Gm)$ denotes the limit set in $\pa \Hb^3$ of a
Kleinian group $\Gm$. 

The limit set of the totally geodesic surface $(\pa
\widetilde{M})_\ell$ is a round circle. Thus the convex hull of $(\pa
\widetilde{M})_\ell$ is a geodesic hyperplane in $\Hb^3$. Since the
convex hull of $\widetilde{\Sigma}_i$ lies to one side of this
hyperplane, we deduce that
\[ \C(\widetilde{\Sigma}_i) \cap (\pa \widetilde{M})_{\ell} = \C \left
   (\Lambda(\pi_1(\widetilde{\Sigma}_i)) \cap \Lambda(\pi_1((\pa
   \widetilde{M})_\ell)) \right ). \] 
By a standard theorem on limit sets (see for example Theorem 3.14
of~\cite{MatsuTani:1998}), the intersection of the limit sets of two
geometrically finite subgroups is the limit set of their intersection,
plus perhaps a set $P$ of parabolic fixed points. By our current
assumption that $M$ has no cusps, $P = \emptyset$. Therefore
\[ \C(\widetilde{\Sigma}_i) \cap (\pa \widetilde{M})_\ell
    = \C \left (\Lambda \left (\pi_1(\widetilde{\Sigma}_i) \cap
    \pi_1((\pa \widetilde{M})_\ell) \right ) \right ). \]
Consider a loop $x$ on $(\pa \widetilde{M})_\ell$ based at $\ast$
which lies in $\pi_1(\widetilde{\Sigma}_i) \cap \pi_1((\pa
\widetilde{M})_\ell)$. Choose a path $p$ from the basepoint $\ast$ on
$(\pa \widetilde{M})_\ell$ to a point on $\widetilde{\Sigma}_i$. Then
take a loop $\sigma$ on $\widetilde{\Sigma}_i$ running through this
point such that $p \sigma p^{-1}$ is homotopic to $x$ relative to
$\ast$. Use this homotopy to define an annulus $f \co A
\looparrowright \widetilde{M}$ with boundary $\sigma \cup x$, obtained
by identifying two sides of the square with sides $x, p, \sigma$ and
$p^{-1}$.

Put $f \co A \looparrowright \widetilde{M}$ into general position with
respect to the embedded surface $\widetilde{\Sigma}_i$, and assume
that the annulus $f(A)$ is transverse to $\widetilde{\Sigma}_i$ in a
neighbourhood of $f(\pa A)$.  It follows that the preimage of $f(A)
\cap \widetilde{\Sigma}_i$ consists of simple closed curves in $A$.
We can remove any inessential simple closed curves $\al$ from $A$
since $\widetilde{\Sigma}_i$ is incompressible: $f(\al)$ will bound
discs in $f(A)$ and therefore in $\widetilde{\Sigma}_i$; these discs
cobound a ball across which $f(A)$ may be homotoped to remove the
intersection.

Therefore we may assume that the preimage of $f(A) \cap
\widetilde{\Sigma}_i$ in $A$ consists of core curves in $A$. Take the
core curve $\al$ nearest the boundary component of $A$ mapping to $x$
in $(\pa \widetilde{M})_\ell$. Let $f \co A_x \looparrowright
\widetilde{M}$ denote the subannulus of $f \co A \looparrowright
\widetilde{M}$ lying between $f(\al)$ and $f(x)$. The interior of this
annulus $f(A_x)$ is disjoint from $\widetilde{\Sigma}_i$, and so it
lies in $\widetilde{M}^\ast_i$.

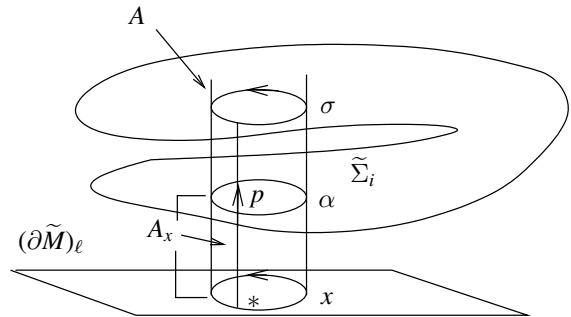
\begin{figure}[ht!]\small
\begin{center}
\begin{picture}(0,0)%
\includegraphics{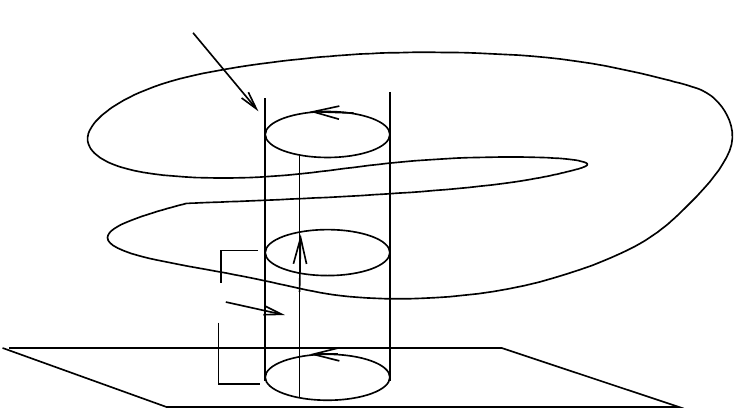}%
\end{picture}%
\setlength{\unitlength}{4144sp}%
\begin{picture}(3361,1862)(1331,-1334)
\put(1375,-925){\makebox(0,0)[lb]{\smash{$(\partial \widetilde{M})_\ell$}}}
\put(3188,-1254){\makebox(0,0)[lb]{\smash{$x$}}}
\put(3173,-676){\makebox(0,0)[lb]{\smash{$\alpha$}}}
\put(3182,-113){\makebox(0,0)[lb]{\smash{$\sigma$}}}
\put(2041,420){\makebox(0,0)[lb]{\smash{$A$}}}
\put(2153,-856){\makebox(0,0)[lb]{\smash{$A_x$}}}
\put(3361,-526){\makebox(0,0)[lb]{\smash{$\widetilde{\Sigma}_i$}}}
\put(2746,-1276){\makebox(0,0)[lb]{\smash{$\ast$}}}
\put(2769,-646){\makebox(0,0)[lb]{\smash{$p$}}}
\end{picture}%
\caption{Homotopy annulus from $\sg$ to $x$} 
\end{center}
\end{figure}

Any such annulus can be homotoped into the characteristic submanifold
$\Psi_i$ of $\widetilde{M}^\ast_i$ by Johannson's Enclosing
Theorem~\cite{Johannson:1979}. By construction, it can be homotoped
into the subset $\Psii \subset \Psi_i$. Thus the curve $x$ on $(\pa
\widetilde{M})_\ell$ can be isotoped into $\Psii \cap (\pa
\widetilde{M})_\ell$. Indeed, since we assumed that the boundary
components of $\Psii \cap \pa \widetilde{M}$ were geodesic, the
geodesic straightening of $x$ already lies inside the intersection
$\Psii \cap (\pa \widetilde{M})_\ell$. We deduce that
$\C(\widetilde{\Sigma}_i) \cap (\pa \widetilde{M})_\ell$ lies inside
$\Psii \cap (\pa \widetilde{M})_\ell$.

Conversely, consider any loop $x$ in $\Psii \cap (\pa
\widetilde{M})_\ell$ based at $\ast$. By definition of $\Psii$, it is
homotopic relative to $\ast$ to a loop of the form $p \sigma p^{-1}$,
where $p$ is some path from $\ast$ to $\widetilde{\Sigma}_i$ and
$\sigma$ is a loop in $\widetilde{\Sigma}_i$. Thus $x \in
\pi_1(\widetilde{\Sigma}_i) \cap \pi_1((\pa \widetilde{M})_\ell)$, and
the geodesic representative of $x$ lies in \[\C(\widetilde{\Sigma}_i)
\cap (\pa \widetilde{M})_\ell = \C \left (\Lambda \left
(\pi_1(\widetilde{\Sigma}_i) \cap \pi_1((\pa \widetilde{M})_\ell)
\right ) \right ).\] Since the boundary curves of both $\Psii \cap
(\pa \widetilde{M})_\ell$ and $\C(\widetilde{\Sigma}_i) \cap (\pa
\widetilde{M})_\ell$ are geodesic, $x$ also lies inside
$\C(\widetilde{\Sigma}_i) \cap (\pa \widetilde{M})_\ell$. This gives
the opposite inclusion. By repeating both directions of this argument
for every component of $\pa \widetilde{M}$, we obtain the required
equality $\Psii \cap \pa \widetilde{M} = \C(\widetilde{\Sigma}_i) \cap
\pa \widetilde{M}$.\endproof

The curves which occur either as cut curves defined in
\fullref{def:cutcurves} or as boundary curves of the set
described in \fullref{prop:hull_equality} will play a special
role in our argument. Sometimes we consider some additional curves as
well.

\begin{Def}
\label{def:charcurves}
We define a \textsl{boundary curve} for $\widetilde{\Sigma}_i$ to be
an essential simple closed geodesic $x$ in $\pa \widetilde{M}$ which
is a component of $\pa (\Psii \cap \pa \widetilde{M})$. We say that an
essential simple closed geodesic in $\pa \widetilde{M}$ is a
\textsl{characteristic curve} for $\widetilde{\Sigma}_i$ if it is a
boundary curve for $\widetilde{\Sigma}_i$ or a cut curve for
$\widetilde{\Sigma}_i$.
\end{Def}

Let $X'_i$ denote the set of characteristic curves for
$\widetilde{\Sigma}_i$. Then $X'_i$ is a possibly empty finite
collection of disjoint essential simple closed geodesics on $\pa
\widetilde{M}$.

\begin{Def}
\label{def:distcurves}
Suppose that $X'_i \cap (\pa \widetilde{M})_\ell$ were empty for some
boundary component $(\pa \widetilde{M})_\ell$ of $\widetilde{M}$. Then
we choose an essential simple closed geodesic $x_\ell \subset (\pa
\widetilde{M})_\ell$ and replace $X'_i$ by $X'_i \cup \{x_\ell\}$. We
call $x_{\ell}$ an \textsl{additional curve}.

Extend each $X'_i$ to a set $X_i$ which is the union of $X'_i$ and an
additional curve $x_\ell$ for every boundary component requiring
one. Then $X_i \cap (\pa \widetilde{M})_\ell$ is non-empty for all $i$
and $\ell$. We call $X_i$ the \textsl{set of distinguished curves} for
$\Psii$. It consists of the characteristic curves for
$\widetilde{\Sigma}_i$ and the additional curves chosen as above. If
$x$ is a curve in $X_i$, we say it is a \textsl{distinguished curve}.
\end{Def}

We will also require the following definition.

\begin{Def}
\label{def:quintersecting}
Let $\al$ be an immersed arc on $(\pa \widetilde{M})_\ell$. We say
that $\al$ is a \textsl{quintersecting arc} if it intersects every
distinguished curve $x\in \bigcup_i X_i$ on $(\pa \widetilde{M})_\ell$
at least five times.
\end{Def}

\section{Ensuring large intersection with $X_i$}

\label{sec:crossing}

In this section we use the definitions above and show that certain
curves whose images bound discs in an attached handlebody have large
geometric intersection number with each distinguished curve in every
set $X_i$.

\begin{Prop}
\label{prop:masurext}
Let $M \neq B^3$ be a simple 3--manifold containing a closed,
connected, orientable, immersed surface $\eta \co \Sigma
\looparrowright M$ with a finite cover which lifts to a non-peripheral
embedded incompressible surface in a degree $d$ cover $\pi \co
\widetilde{M} \to M$. Let $H$ be a handlebody. Suppose that $h \co \pa
H \to \pa H$ is a pseudo-Anosov map whose stable lamination is of full
type and that $\phi' \co (\pa M)_\ell \to \pa H$ is a homeomorphism
from some boundary component $(\pa M)_\ell$ of $M$ to $\pa H$.

Then, for all $n \in \Nb$, there exists $N_{\textrm{min}}(h,\phi') \in
\Nb$ depending on the maps $\eta$, $\pi$, $h$ and $\phi'$ and on the
integers $n$ and $d$ such that every attaching homeomorphism $\phi =
h^{N} \circ \phi' \co (\pa M)_\ell \to \pa H$ with $N \geq
N_{\textrm{min}}(h,\phi')$ has the following property.

$[\ast]$ Any connected essential (not necessarily embedded) curve
$\rho \subset (\pa \widetilde{M})_\ell$ such that $\phi \circ
\pi(\rho)$ bounds a disc in $H$ can be subdivided into at least $n$
subarcs, each of which is a quintersecting arc.
\end{Prop}

\proof Take the component of $\pa \widetilde{M}$ containing
$\rho$ and denote it by $(\pa \widetilde{M})_\ell$. Thus $(\pa
\widetilde{M})_\ell$ projects to $(\pa M)_\ell$. The homeomorphism
$\phi' \co (\pa M)_\ell \ra \pa H$ allows us to define a cover
$\pi_{\phi'}$ of $\pa H$ by
\[ \pi_{\phi'} := \phi' \circ \pi \co \pa \widetilde{M}_{\ell} \ra \pa
H \] We write $\widetilde{\pa H}=\pa \widetilde{M}_\ell$. Note that in
 general $\pi_{\phi'}$ does not extend to a cover of the handlebody
 $H$. We then want to be able to lift a homeomorphism $h \co \pa H \to
 \pa H$ to a homeomorphism $\widetilde{h} \co \widetilde{\pa H} \to
 \widetilde{\pa H}$. The necessary and sufficient condition for this
 to be possible is that the map $h_\ast$ induced by $h$ at the level
 of the fundamental group should satisfy $h_\ast \Gm=\Gm$, where $\Gm
 \leq \pi_1(\pa H)$ is the subgroup
\[ \Gm = (\pi_{\phi'})_\ast \left ( \pi_1(\widetilde{\pa H}) \right
). \] 
This subgroup is of finite index equal to $\bar{d} \leq d$, where $d$
is the degree of the cover $\pi \co \widetilde{M} \to M$. The group
$\pi_1(\pa H)$ is finitely generated, hence it has finitely many
subgroups of index $\bar{d}$, and the automorphism $h_\ast$ must
permute them. Such a permutation has order $\nu$ dividing $t!$ where
\[ t = \max_{\bar{d}\leq d}{\left \{ \mbox{number of subgroups of
    $\pi_1(\pa H)$ of index $\bar{d}$} \right \}} < \infty
\]
Thus for any $h$, the induced map $(h^{t!})_\ast = (h_\ast)^{t!}$
preserves $\Gm$. This implies that $h^{t!}$ lifts to a homeomorphism
\[\widetilde{h^{t!}} \co \widetilde{\pa H} \ra \widetilde{\pa H}.\]
Consider the map $\phi = h^N \circ \phi'$, and take $n \in \Nb$. We
may write $N$ as $N't! + r$, where $0 \leq r < t!$. Write $\phi$ as
\[ \phi = (h^{t!})^{N'} \circ (h^r \circ \phi'),\]
noting that none of $N', t$ or $r$ depends on $\phi'$.

Suppose that $\psi \co (\pa M)_\ell \ra \pa H$ ranges over the
homeomorphisms in the set
\[ \Phi_{t!} = \{\phi', h\phi', h^2 \phi, \ldots, h^{t!-1}\phi'\}.\]
If for every $\psi \in \Phi_{t!}$ we can obtain an integer
$N_{\textrm{min}}(h^{t!},\psi)$ which satisfies the requirements of
the proposition, then we can take
\[ N_{\textrm{min}}(h,\phi') = t! \times \left ( \max_{\psi \in
  \Phi_{t!}} N_{\textrm{min}}(h^{t!},\psi) \quad + \quad 1 \right ).\]
Indeed, if we take $N \geq N_{\textrm{min}}(h,\phi')$ with
$N_{\textrm{min}}(h,\phi')$ defined as above, then for
\[ h^N \circ \phi' = (h^{t!})^{N'} \circ (h^r \circ \phi'), \]
we have 
\[ N' \quad \geq \quad \frac{N}{t!} - 1 \quad \geq \quad
N_{\textrm{min}}(h^{t!},\psi)\] with $\psi=h^r \circ \phi'$, so that
the conclusion holds.

Thus from now on we will assume that the map $h$ lifts to a
pseudo-Anosov homeomorphism $\widetilde{h} \co \widetilde{\pa H} \to
\widetilde{\pa H}$. The stable lamination $\widetilde{\lambda}^S$ of
$\widetilde{h}$ is $\pi_{\phi'}^{-1}(\lambda^S)$. 

We define $\widetilde{\phi} \co \pa \widetilde{M}_\ell \ra
\widetilde{\pa H}$ to be $\widetilde{h}^N$ using $\widetilde{\pa H} =
\pa \widetilde{M}_\ell$. Then
\[\pi_{\phi'} \circ \widetilde{\phi} = \pi_{\phi'} \circ
 \widetilde{h}^N = h^N \circ \pi_{\phi'} = h^N \circ \phi' \circ \pi 
 = \phi \circ \pi \]
so that gluing commutes with the projection map.

By our choice of map $h$, $\lambda^S \in P\Wc(H)^c$. Thus the
lamination $\lambda^S$ satisfies all strict triangle inequalities on
each pants component $P$ of some handlebody pants decomposition
$\Pc=\Pc(\lambda^S)$. As $\lambda^S$ has no self-intersections and
never spirals into a closed leaf, the strict triangle inequalities
imply that $\lambda^S$ has no waves and has all possible seams.
Moreover, the union $\lambda^S \cup \Pc$ splits $\pa H$ up into
simply-connected regions, either square or hexagonal.

\begin{figure}[ht!]\small
\begin{center}
\begin{picture}(0,0)%
\includegraphics{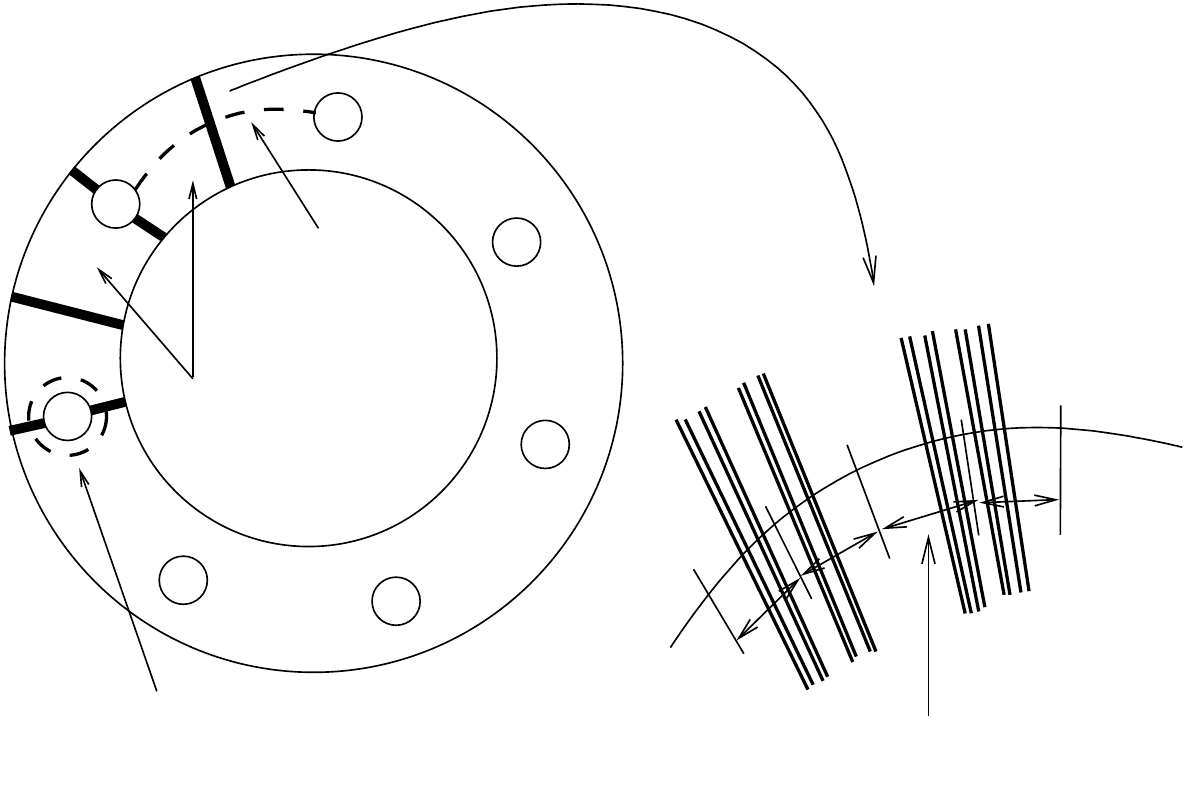}%
\end{picture}%
\setlength{\unitlength}{4144sp}%
\begin{picture}(5418,3598)(136,-3211)
\put(1253,-953){\makebox(0,0)[lb]{\smash{in $\phi \circ \pi(\rho)$}}}
\put(1254,-796){\makebox(0,0)[lb]{\smash{lift of a wave}}}
\put(1951,-391){\makebox(0,0)[lb]{\smash{lift of $P$}}}
\put(3871,142){\makebox(0,0)[lb]{\smash{close-up of intersection}}}
\put(923,-1463){\makebox(0,0)[lb]{\smash{hexagonal regions}}}
\put(2897,-2491){\makebox(0,0)[lb]{\smash{$\widetilde{\phi}(\rho)$}}}
\put(3481,-3053){\makebox(0,0)[lb]{\smash{subarcs containing intersections}}}
\put(3781,-1591){\makebox(0,0)[lb]{\smash{$\widetilde{\lambda}^S$}}}
\put(496,-2986){\makebox(0,0)[lb]{\smash{lift of component of $\phi \circ \pi(\rho)$}}}
\put(496,-3166){\makebox(0,0)[lb]{\smash{which is parallel to some 
$\partial P$}}}
\end{picture}%
\caption{Every lift of a wave intersects 
         $\widetilde{\lambda}^S$ many times}
\label{fig:lifting_to_sc_regions}
\end{center}
\end{figure}
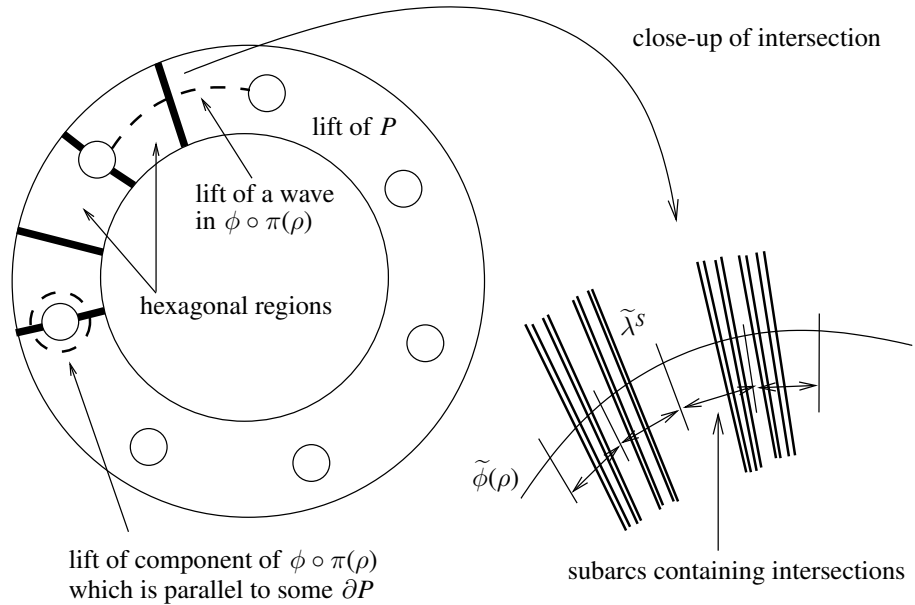

In the cover $\widetilde{\pa H}$, these regions must lift to
simply-connected regions. \fullref{fig:lifting_to_sc_regions}
shows an example of a possible lift of some $P \in \Pc$.

We now carry out the technical part of the proof. Take a curve $\rho$
in $\pa \widetilde{M}$ (connected but not necessarily simple) such
that $\phi \circ \pi(\rho)= \pi_{\phi'} \circ \widetilde{\phi}(\rho)$
bounds a disc in $H$. This disc is proper and essential in $(H,\pa H)$
but need not be embedded. By \fullref{lem:bdhwave}, either $\phi
\circ \pi(\rho)$ is freely homotopic into $P \in \Pc$ and hence
separates two boundary components of the \Hpants $P$, or it has an
immersed wave with respect to $\Pc$. In either case, $\phi \circ
\pi(\rho)$ intersects $\lambda^S$.

This intersection persists when we lift to $\widetilde{\pa H}$:
$\widetilde{\phi}(\rho)$ must intersect
$\widetilde{\lambda}^S$. Indeed, the lift of the immersed wave (or
boundary-parallel curve) starts in one of the square or hexagonal
regions and must leave by one of the edges lying in
$\widetilde{\lambda}^S$. Furthermore, since the lamination
$\widetilde{\lambda}^S$ has local geometry modelled on a Cantor set
$\times\mbox{ }\Rb$, the lift of the wave intersects
$\widetilde{\lambda}^S$ uncountably many times near this point.

Apply the map $\widetilde{\phi} = \widetilde{h}^N$ to a simple closed
curve $x$ in $(\pa \widetilde{M})_\ell$. Denoting the geodesic
straightening of a curve's image by $(-)^\ast$, we have
\[ \left (\widetilde{\phi} (x)\right )^\ast \quad = \quad
(\widetilde{h}^N (x))^\ast \quad \subset \quad \widetilde{\pa H}.\] As
$N \ra \infty$, $(\widetilde{h}^N(x))^\ast \ra \widetilde{\lambda}^S$
in the Gromov--Hausdorff sense (e.g.~\cite{Kapovich:2001}, $\S$11.15).

Write $W(x)$ for the smallest integer (depending on $\eta$, $\pi$, $h$
and $\phi'$) such that whenever $N \geq W(x)$, the following property
holds. This integer will exist by the convergence above. 

\begin{Property} For every lift $\widetilde{\xi}$ of a seam $\xi$ in
  an \Hpants $P$ of $\Pc$ such that $\widetilde{\lambda}^S$ has a
subarc homotopy equivalent to $\widetilde{\xi}$ keeping its endpoints
on $\pa P$, the curve $(\widetilde{h}^N(x))^\ast$ has a subarc
homotopy equivalent to $\widetilde{\xi}$ keeping its endpoints on $\pa
P$.
\end{Property}

Note that there may not be a uniform bound $W$ such that, for all
simple closed curves $x$, $W(x) \leq W$. Consider however the set of
distinguished curves $X_i$ in $\pa \widetilde{M}$ from
\fullref{def:distcurves}. Since there are finitely many curves
in each $X_i$, and finitely many lifts $\widetilde{\Sigma}_i$ from
which we obtain a set of distinguished curves $X_i$, there are
finitely many curves in all of these sets put together.

Write $X_{(\pa \widetilde{M})_\ell}$ for $\cup_i X_i \cap (\pa
\widetilde{M})_\ell$.  Since $X_{(\pa \widetilde{M})_\ell}$ is finite,
we can define $W_{X, \ell}$ to be the maximum of the integers $W(x)$
such that $x \in X_{(\pa \widetilde{M})_\ell}$. Then, for $N \geq
W_{X, \ell}$ and a geodesic $x \in X_{(\pa \widetilde{M})_\ell}$, the
curves
\[(\widetilde{\phi}(\rho))^\ast \equiv (\widetilde{h}^N(\rho))^\ast
\qquad \mbox{and} \qquad (\widetilde{h}^N(x))^\ast\] have non-zero
geometric intersection number. Thus $\rho$ and $x$ intersect, proving
a weak form of the proposition.

More generally, take a point $p$ in the intersection of
$\widetilde{\lambda}^S$ with
$\widetilde{\phi}(\rho)=\widetilde{h}^N(\rho)$. As $N$ tends to
infinity, the intersection of $(\widetilde{\phi}(x))^\ast =
(\widetilde{h}^{N}(x))^\ast$ with $\widetilde{h}^N(\rho)$ converges to
$\widetilde{\lambda}^S \cap \mbox{ }\widetilde{h}^N(\rho)$ around $p$
in the sense of Gromov--Hausdorff. In particular, the intersections
tend to the Cantor set limit.

To complete the proof, take $n \in \Nb$. Any curve
$\widetilde{\phi}(\rho)$ as above has at least $n$ subarcs each of
which intersects $\widetilde{\lambda}^S$ at least five times, possibly
all on the same lift of a wave. This implies that for each
distinguished curve $x$ there exists $N_{\textrm{min}}(x)$ such that
$\widetilde{\phi}(\rho)$ has at least $n$ subarcs each of which
intersects $\phi(x)$ at least five times. 

Since there are finitely many distinguished curves and their images
all converge to $\widetilde{\lambda}^S$, we deduce that there exists
$N_{\textrm{min}}$ such that $\widetilde{\phi}(\rho)$ has $n$ subarcs
each of which intersects all distinguished curves at least five
times. For example, pick $5n$ consecutive intersections between
$\widetilde{\phi}(\rho)$ and $\widetilde{\lambda}^S$, a number $\ep>0$
small enough so that the distance between any two such intersections
is greater than $2\ep$, and take $N_{\textrm{min}}$ large enough so
that for every distinguished curve $x$, the geodesic $(\phi(x))^\ast$
has an intersection with $\widetilde{\phi}(\rho)$ within an
$\ep$--neighbourhood of each intersection with
$\widetilde{\lambda}^S$. Thus the same subarcs may be chosen so that
they work for all distinguished curves.

Moreover, the integer $N_{\textrm{min}}$ can be taken to be
independent of $\rho$ since there are only finitely many types of
square or hexagonal regions to consider. We conclude the proof by
noting that this property is preserved under
$\widetilde{\phi}^{-1}$. Thus there exists $N_{\textrm{min}}$ so that
every $\rho$ can be subdivided into at least $n$ subarcs, each of
which is quintersecting.\endproof

\section{Gathering area outside $\C(\widetilde{\Sigma}_i)$}
\label{sec:area_outside_hull}

In Sections~\ref{sec:cutlocus} and~\ref{sec:crossing} we obtained two
useful results. \fullref{prop:hull_equality} expresses the
intersection of the convex hull of a surface $\widetilde{\Sigma}_i$
with $\pa \widetilde{M}$ as the intersection of a certain collection
of $I$--bundles $\Psii$ with $\pa \widetilde{M}$.
\fullref{prop:masurext} gives conditions guaranteeing that
certain curves in $\pa \widetilde{M}$ have many subarcs, each
intersecting every curve in the set $X_i$ of distinguished curves for
$\Psii$.

In this section and \fullref{sec:incomp} we combine these results
to show that all but a small number of these subarcs containing
enough of these intersections will contribute a fixed quantum of area
to a lift of a planar surface lying in $M$. This surface will be of
the following kind.

Let $f \co Q \looparrowright M$ be a least-area planar surface
immersed in $M$, with $k \geq 2$ geodesic boundary components
$f(q_1),\ldots,f(q_k)$ lying on $\pa M$ and a geodesic boundary
component $f(q_0)$ lying in $\C(\Sigma)$ and homotopic to a curve in
$\Sigma$. This surface $Q$ inherits a negatively curved metric from
$M$ (see Theorem 5.5, $\S$V.A,~\cite{GHL:1987}) of curvature at most
$-1$. Assume also that $Q$ is homotopically incompressible and
homotopically $\pa$--incompressible.

Basmajian~\cite{Basmajian:1994} showed that the totally geodesic
boundary $\pa M$ of a hyperbolic manifold such as $M$ has a collar of
width some fixed constant $U$ depending only on $\{\chi((\pa
M)_\ell)\}$. This leads us to the following definitions of a \Ucollar
and of minimal \Uthin arcs in a planar surface immersed in $M$.

\begin{Def}
\label{def:ucollar}
For $M$, $\Sigma$ and $Q$ as above, consider a curve $q_j$ ($j \neq
0$). The \textsl{U--collar} of $q_j$ consists of all points $\zeta \in
Q$ such that the distance $d(\zeta,q_j)$ is at most $U$ in the path
metric on $Q$ obtained from the metric on $M$.
\end{Def}

By Basmajian's result these collars are disjoint for $q_j$ ($j \neq
0$). For if not, we can find an embedded essential arc in $Q$ running
between curves $q_{j_1}$ and $q_{j_2}$ which lies in the \Ucollar of
$\pa M$ and so can be homotoped into $\pa M$. This contradicts
homotopic $\pa$--incompressibility of $Q$. However, the curve $f(q_0)$
may approach some $f(q_j)$ at a distance less than $U$, and so we
define $\de$--thin arcs for $Q$.

\begin{Def}
\label{def:thinpart}
Define a \textsl{minimal U--thin arc} in $f \co Q \looparrowright M$
to be a geodesic arc $\xi$ properly embedded in $Q$ with endpoints on
$q_0$ and some $q_j$ ($j \neq 0$), such that (i) the length of $\xi$
is at most $U$ in the path metric on $Q$; (ii) $\xi$ is of minimal
length among embedded arcs homotopy equivalent to it keeping their
endpoints on the boundary.
\end{Def}

\begin{figure}[ht!]\small
\begin{center}
\begin{picture}(0,0)%
\includegraphics{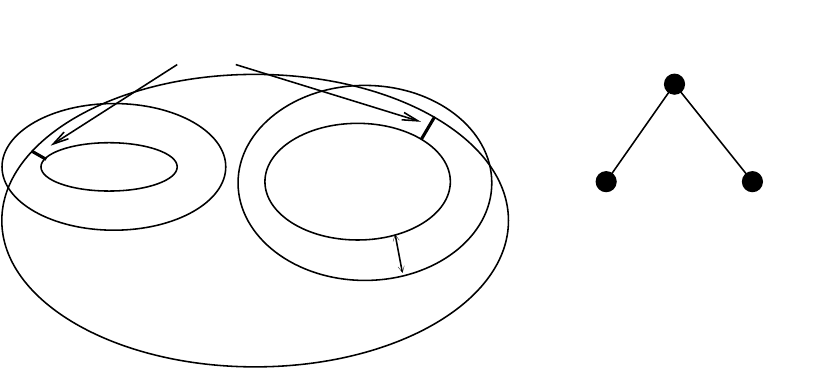}%
\end{picture}%
\setlength{\unitlength}{3947sp}%
\begin{picture}(3935,1754)(665,-1218)
\put(2616,-696){\makebox(0,0)[lb]{\smash{$U$}}}
\put(1470,299){\makebox(0,0)[lb]{\smash{$U$\negthinspace--thin arc}}}
\put(1231,-271){\makebox(0,0)[lb]{\smash{$q_2$}}}
\put(2229,-196){\makebox(0,0)[lb]{\smash{$q_1$}}}
\put(1186,-872){\makebox(0,0)[lb]{\smash{$Q$}}}
\put(1449,-1119){\makebox(0,0)[lb]{\smash{$q_0$}}}
\put(3504,-527){\makebox(0,0)[lb]{\smash{$v_2$}}}
\put(4216,-526){\makebox(0,0)[lb]{\smash{$v_1$}}}
\put(3857,238){\makebox(0,0)[lb]{\smash{$v_0$}}}
\put(3587,404){\makebox(0,0)[lb]{\smash{$\mathcal{G}_Q$}}}
\end{picture}%
\caption{Definition of minimal \Uthin arcs and example of a graph}
\end{center}
\end{figure}
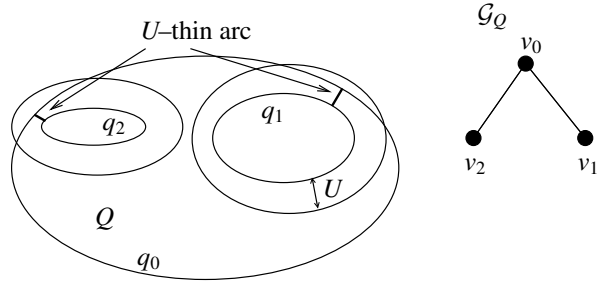

Note that these arcs are disjoint. For any constant $0 <\de \leq U$,
we may similarly define \textsl{minimal $\mathbf{\de}$--thin arcs} for
$f \co Q \looparrowright M$. Moreover, if $R$ is a cover of $Q$ which
lifts to some finite cover $\widetilde{M} \to M$, the set of minimal
$\de$--thin arcs of $R$ consists of all arcs in $R$ which project to a
minimal $\de$--thin arc of $Q$.

\begin{Lem}
\label{lem:thinpart}
Let $\tilde{f} \co R \looparrowright \widetilde{M}$ be a (minimal)
surface immersed in $\widetilde{M}$ such that its projection
$Q=\pi(R)$ is a least-area planar surface in $M$ as described above,
with $k+1 \geq 3$ boundary components. Suppose that $d$ is the degree
of the cover $\pi \co \widetilde{M} \to M$. Then, for any $0<\de \leq
U$, there are at most $2(k-1)d$ minimal $\de$--thin arcs in $R$.
\end{Lem}

\proof By construction it is enough to show that $Q$ has at
most $2(k-1)$ minimal $\de$--thin arcs, where $k \geq 2$. Consider
therefore the graph $\Gc_Q$ obtained from $Q$ by taking a vertex for
every boundary component of $Q$ and an edge for every minimal
$\de$--thin arc. Then $\Gc_Q$ is planar, and we claim it has at most
$2(k-1)$ edges.

We prove this by induction on $k$. If $k=2$, there are at most two
edges, each of which runs from the vertex $v_0$ representing $q_0$ to
one of the other vertices $v_1$ or $v_2$. For, since each $f(q_j)$ is
a geodesic curve in the negatively-curved surface $Q$, there cannot
be two disjoint adjacent minimal $\de$--thin arcs without some
non--simply-connected region of $Q$ separating them.

If $k>2$, we can always find a vertex with at most one edge connecting
it to $v_0$. For, suppose $v$ is some vertex with more than such one
edge, and let $\xi_1$ and $\xi_2$ be minimal $\de$--thin arcs
represented by two of these edges.

Consider cutting $Q$ along $\xi_1 \cup \xi_2$. Since $\xi_1$ and
$\xi_2$ have some non-simply connected region of $Q$ lying between
them, both components of $Q-(\xi_1 \cup \xi_2)$ are planar surfaces
of Euler characteristic at most zero. Hence each such surface has
Euler characteristic whose absolute value is strictly less than
$|\chi(Q)|$.

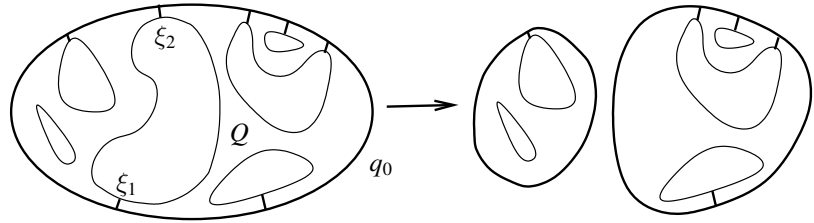
\begin{figure}[ht!]\small
\begin{center}
\begin{picture}(0,0)%
\includegraphics{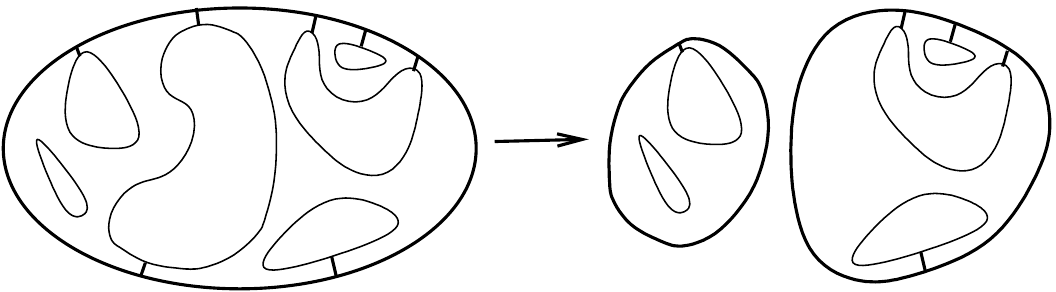}%
\end{picture}%
\setlength{\unitlength}{3947sp}%
\begin{picture}(5060,1384)(90,-628)
\put(751,-437){\makebox(0,0)[lb]{\smash{$\xi_1$}}}
\put(2344,-311){\makebox(0,0)[lb]{\smash{$q_0$}}}
\put(1466,-124){\makebox(0,0)[lb]{\smash{$Q$}}}
\put(993,485){\makebox(0,0)[lb]{\smash{$\xi_2$}}}
\end{picture}%
\caption{Cutting $Q$ along some $\xi_1$ and $\xi_2$}
\end{center}
\end{figure}

Therefore we may argue by induction on $|\chi(Q)|$, choosing one of
these new planar surfaces each time. But $|\chi(Q)|$ is finite so this
process must terminate. We eventually find a boundary component $q_j$
of $Q$ incident with at most one minimal $\de$--thin arc;
equivalently, we find a vertex $v_j$ with at most one edge connecting
it to $v_0$.

Removing this boundary component causes at most one pair of the
minimal $\de$--thin arcs to become adjacent. Since the curvature of
$Q$ is negative, these will combine into a single minimal arc.

In the graph, this corresponds to removing the vertex with at most one
edge, and amalgamating at most one pair of edges into a single
edge. In this way we obtain a graph with one fewer vertex and at most
two fewer edges. By induction, $\Gc_Q$ has at most $2(k-1)$ edges.\endproof

We have now shown that, for any $0<\de \leq U$, the surface $R$ has at
most $2(k-1)d$ minimal $\de$--thin arcs. These will provide the
possible exceptional cases in the argument which follows in
\fullref{sec:incomp}. Before this, we will determine a couple of
constants which will appear in future formulae.

\begin{Def}
For $\Sigma$, $M$ and $\widetilde{M}$ as described at the beginning of
\fullref{sec:cutlocus}, recall the set $X_i$ of distinguished
curves from \fullref{def:distcurves}.

Consider all geodesic curves $x_{i,j} \in X=\bigcup_i X_i$. They each
have an annular collar in $\pa \widetilde{M}$ of some strictly
positive width $w_{i,j}$. For fixed $i$ these collars are assumed to
be disjoint. Define a constant $\ep>0$ by $\ep= \frac{1}{4} \min_{i,j}
w_{i,j}$. We say that $\ep$ is the \textsl{minimum disjoint annular
width} for $X$.
\end{Def}

It is useful to consider a special case. Choose a boundary component
$(\pa \widetilde{M})_\ell$ and suppose that each set $X_i$ has some
boundary curves on $(\pa \widetilde{M})_\ell$ in it. By
Definitions~\ref{def:charcurves} and~\ref{def:distcurves} this is
equivalent to assuming that, for each $i$, the set
$\C(\widetilde{\Sigma}_i) \cap (\pa \widetilde{M})_\ell$ is a
non-empty proper subset of $(\pa \widetilde{M})_\ell$.


\begin{figure}[ht!]\small
\begin{center}
\begin{picture}(0,0)%
\includegraphics{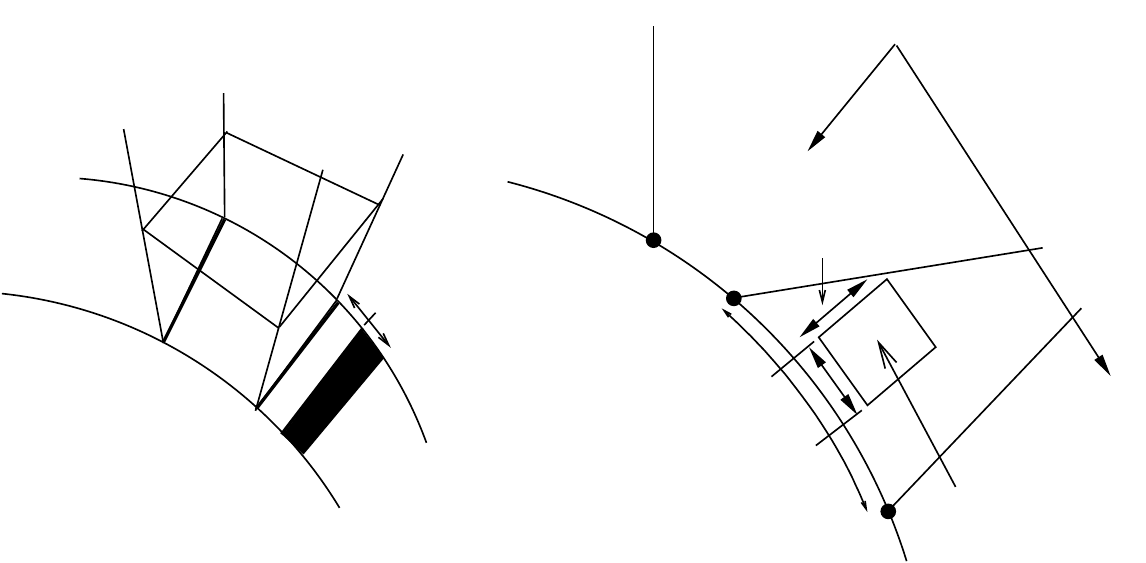}%
\end{picture}%
\setlength{\unitlength}{4144sp}%
\begin{picture}(5231,2570)(86,-1757)
\put(1648,-1229){\makebox(0,0)[lb]{\smash{$A_{i,2}$}}}
\put(3457,-873){\makebox(0,0)[lb]{\smash{$\epsilon$}}}
\put(3896,681){\makebox(0,0)[lb]{\smash{convex hull of $\widetilde{\Sigma}_i$}}}
\put(1185,221){\makebox(0,0)[lb]{\smash{convex hull}}}
\put(4456,-1501){\makebox(0,0)[lb]{\smash{cross-section }}}
\put(4456,-1681){\makebox(0,0)[lb]{\smash{of collar}}}
\put(205,-381){\makebox(0,0)[lb]{\smash{$\partial{\widetilde{M}}$}}}
\put(702,-895){\makebox(0,0)[lb]{\smash{$x_{i,1}$}}}
\put(1090,-1190){\makebox(0,0)[lb]{\smash{$x_{i,2}$}}}
\put(3661,-333){\makebox(0,0)[lb]{\smash{$\geq \delta$}}}
\put(2953,-454){\makebox(0,0)[lb]{\smash{$x_{i,1}$}}}
\put(3917,-1665){\makebox(0,0)[lb]{\smash{$x_{i,3}$}}}
\put(2701,-61){\makebox(0,0)[lb]{\smash{$\partial \widetilde{M}$}}}
\put(3229,-722){\makebox(0,0)[lb]{\smash{$x_{i,2}$}}}
\put(1732,-570){\makebox(0,0)[lb]{\smash{$\epsilon$}}}
\put(1847,-690){\makebox(0,0)[lb]{\smash{$\epsilon$}}}
\put(3655,-1116){\makebox(0,0)[lb]{\smash{$\epsilon$}}}
\end{picture}%
\caption{Obtaining a collar disjoint from the convex hull}
\end{center}
\end{figure}
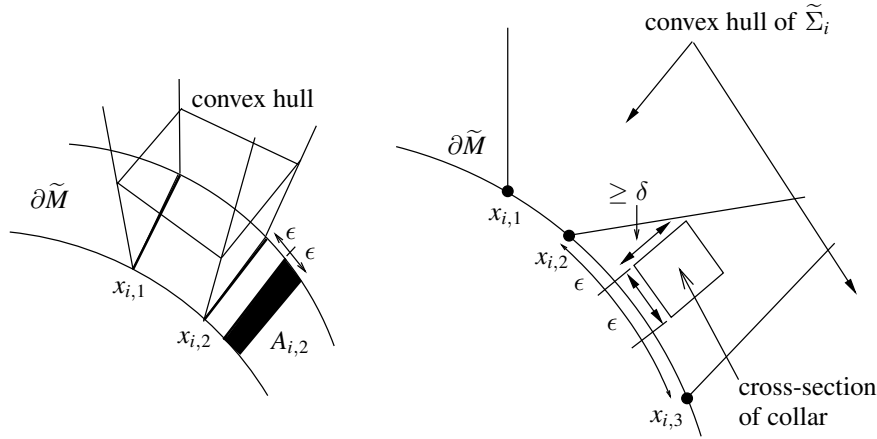

Choose a lift $\widetilde{\Sigma}_i$ of $\Sigma$ and consider its
convex hull $\C(\widetilde{\Sigma}_i)$. This intersects $(\pa
\widetilde{M})_\ell$ in subsurfaces of $\widetilde{M}$ whose boundary
curves lie at least a distance $4\ep$ apart. For each boundary curve
$x_{i,j}$, consider an annulus $A_{i,j}$ in $(\pa \widetilde{M})_\ell$
which lies parallel to $x_{i,j}$ and which is disjoint from the convex
hull $\C(\widetilde{\Sigma}_i)$. We assume that the distance between
$x_{i,j}$ and the annulus is $\ep$, and that the annulus has width
$\ep$.

\begin{Lem}
\label{lem:delta}
There exists a constant $\de>0$ such that, for all $i,j$, the annulus
$A_{i,j}$ has a solid collar $A_{i,j} \times I$ of width at least
$\de$, where this collar is disjoint from the convex hull of
$\widetilde{\Sigma}_i$.
\end{Lem}

\proof By taking $\de \leq U$, we know that each $A_{i,j}$ has
a collar of width at least $\de$ in $\widetilde{M}$. Each annulus
$A_{i,j}$ lies outside the appropriate convex hull
$\C(\widetilde{\Sigma}_i) \cap (\pa \widetilde{M})_\ell$, and there
are finitely many such annuli. By compactness for each $i$ and $\ell$,
the annuli $A_{i,j}$ lie at least some distance $2\de_i$ from the
appropriate convex hull. Thus they each have a collar of width $\de =
\min_i \de_i$ which is disjoint from it.\endproof

Take a $\de$--collar of the boundary of $\widetilde{M}$, and denote it
by $N_\de(\pa \widetilde{M})$. By decreasing $\de>0$ if necessary, the
intersection of the convex hull of each $\widetilde{\Sigma}_i$ with
$N_\de(\pa \widetilde{M})$ is a collar on the intersection of that
convex hull with the boundary of $\widetilde{M}$.

We may do this as explained before
\fullref{prop:hull_equality} by making an extremely small
modification (of width $<\de/1000$, say) to extend the convex hull to
a submanifold with smooth boundary. This condition will ensure that
arguments which work for boundary curves also work for any additional
curves. We will deal with cut curves separately later.

\begin{Def}
\label{def:collarwidth}
We say that any constant $\de>0$ satisfying \fullref{lem:delta} and
the above condition on the collar of every convex hull is a
\textsl{minimum collar width} for $X$. The value of this constant
depends only on $M$, $\widetilde{M}$ and $\Sigma$.
\end{Def}

We have now defined two constants depending on $M$, $\widetilde{M}$
and $\Sigma$: the minimum disjoint annular width $\ep$ and the minimum
collar width $\de$. These will enable us to prove our main theorem in
\fullref{sec:incomp}. The idea is to show that a planar surface
intersects many of these annular collars, and hence picks up at least
$\ep \mbox{ sinh }\de$ area each time, leading to a contradiction with
the Gauss--Bonnet theorem.

\section{Proof of main theorem: non-cusped case}
\label{sec:incomp}

We now use the results obtained in Sections~\ref{sec:projwavelike}
to~\ref{sec:area_outside_hull} to deduce our main result in the case
where the manifold $M$ has no cusps. We consider the general case in
\fullref{sec:many_components}. 

\begin{Thm}
\label{thm:sfsubgp}
Suppose that $M$ is a simple 3--manifold with $m\geq 1$ boundary
components $(\pa M)_1,\ldots,(\pa M)_m$ each of genus $\geq
2$. Consider a collection of handlebodies $\Hc=\{H_1,\ldots,H_m\}$
whose genera match those of $\pa M$. Let $M \cup_\phi \Hc$ denote the
closed 3--manifold obtained by gluing each boundary component $(\pa
M)_{\ell}$ to $\pa H_{\ell}$ via a homeomorphism $\phi_{\ell}$.

Suppose moreover that $h_\ell \co \pa H_\ell \ra \pa H_\ell$ is a
pseudo-Anosov homeomorphism with stable lamination of full
type. Given homeomorphisms $\phi'_\ell \co (\pa M)_\ell \to \pa
H_\ell$, there exist integers $(N_\ell)_{\textrm{min}}$ such that if
$\phi_\ell = h_\ell^{N_\ell} \circ \phi'_\ell$ with $N_\ell \geq
(N_\ell)_{\textrm{min}}$ for all $\ell$, the group $\pi_1(M
\cup_{\phi} \Hc)$ contains a surface subgroup.
\end{Thm}

Recall that $M \neq B^3$ is simple if it is compact, irreducible,
$\pa$--irreducible, atoroidal and acylindrical, so $\pa M$ is
incompressible. A surface subgroup is a subgroup which is isomorphic
to the fundamental group of a closed surface of genus at least two. We
defined laminations of full type in \fullref{def:fulltype}.

\proof Any 3--manifold $M \neq B^3$ which is irreducible,
atoroidal and acylindrical and which has non-empty incompressible
boundary is Haken. As a consequence of the proof of Thurston's
geometrization of Haken manifolds, $M$ can be given a hyperbolic
structure with totally geodesic boundary in the sense of
\fullref{sec:intro}.

By Cooper--Long--Reid~\cite{CooperLongReid:1997}, there exists a
$\pi_1$--injective surface $\eta \co \Sigma \looparrowright M$. Thus
$M$ has a surface subgroup. This surface $\Sigma$ has a finite cover
which lifts to a non-peripheral embedded incompressible surface in a
finite cover $\pi \co \widetilde{M} \to M$, which we can further
assume to be a regular finite cover.

It is therefore enough to show the following theorem. We will make
free use of any terminology defined in the preceding sections. We
usually abuse notation by writing $\Hc$ for the union of the
handlebodies $H_1,\ldots,H_m$.

\begin{Thm}
\label{thm:sfs_incomp}
Let $M$ be a hyperbolic 3--manifold with non-empty totally geodesic
boundary. Let $\eta \co \Sigma \looparrowright M$ be any closed,
connected orientable surface with a finite cover which lifts to a
non-peripheral embedded incompressible surface in a finite regular
cover $\pi \co \widetilde{M} \to M$. Let $M \cup_\phi \Hc$ denote the
closed 3--manifold obtained by gluing each boundary component $(\pa
M)_{\ell}$ to $\pa H_{\ell}$ via a homeomorphism $\phi_{\ell}$, where
each $H_\ell$ is a handlebody as in \fullref{thm:sfsubgp}.

Suppose moreover that $h_\ell \co \pa H_\ell \ra \pa H_\ell$ is a
pseudo-Anosov homeomorphism with stable lamination of full type. Given
homeomorphisms $\phi'_\ell \co (\pa M)_\ell \to \pa H_\ell$, there
exist integers $(N_\ell)_{\textrm{min}}$ such that if $\phi_\ell =
h_\ell^{N_\ell} \circ \phi'_\ell$ with $N_\ell \geq
(N_\ell)_{\textrm{min}}$ for all $\ell$, the surface $\Sigma$ stays
$\pi_1$--injective in the resulting manifold $M \cup_\phi \Hc$.

In particular, $\pi_1(\Sigma) \leq \pi_1(M \cup_\phi \Hc)$ is a
surface subgroup.
\end{Thm}

\proof We argue by contradiction. If the surface $\Sigma$ does
not remain $\pi_1$--injective in $M \cup_\phi \Hc$, there is some
essential loop $\Lc$ in $\eta(\Sigma)$ bounding an immersed disc $f
\co D \looparrowright M \cup_\phi \Hc$. By the hypotheses, we can take
each $(N_\ell)_{\textrm{min}}$ large enough to apply a result of
Lackenby~\cite{Lackenby:2001}, and deduce that $\pi_1(\Hc)$ injects
into $\pi_1(M \cup_\phi \Hc)$. In this case, the intersections of the
disc $f(D)$ with $M$ and $\Hc$ can be assumed to have the following
structure.

\begin{Lem}
\label{lem:discstructure}
Suppose $\pi_1(\Hc)$ injects into $\pi_1(M \cup_\phi \Hc)$. For any
loop $\Lc$ in $M$ homotopically trivial in $M \cup_\phi \Hc$, take a
spanning disc $f \co D \looparrowright M \cup_\phi \Hc$ with $\Lc$
running once round $f(D)$. Assume moreover that $D \cap f^{-1}(\Hc)$
consists of a collection of discs which are disjoint in $D$, and that
the number of such discs is minimal. Then $f(D) \cap M$ is a
homotopically incompressible and homotopically $\pa$--incompressible
surface in $M$.
\end{Lem}

Note that here homotopically incompressible means that no
homotopically non-trivial simple closed curve in $D \cap f^{-1}(M)$
maps to a homotopically trivial curve in $M$; homotopically
$\pa$--incompressible means that there is no properly embedded
essential arc in $D \cap f^{-1}(M)$ which maps to an arc which can be
homotoped into $\pa M$ while keeping its endpoints fixed.

\proof We can always find $f \co D \looparrowright M \cup_\phi
\Hc$ such that $f(D) \cap \Hc$ is a disjoint union of discs by
considering the handlebodies in $\Hc$ to be small neighbourhoods of
their core graphs. Indeed, these graphs can be made to intersect
$f(D)$ transversely in isolated points. Thus we can assume that $D
\cap f^{-1}(\Hc)$ is a collection of small discs embedded in $D$, and
that the number of such discs is minimal.

Suppose $f(D) \cap M$ were homotopically compressible in $M$. Then
there would exist a homotopically non-trivial simple closed curve $c$
in $D \cap f^{-1}(M)$ which mapped to a homotopically trivial curve in
$M$. We could then modify $f$ by replacing the interior of $c$ by the
disc it spans in $M$, reducing $|D \cap f^{-1}(\Hc)|$. Since this
number was assumed minimal, we have a contradiction. Thus we see that
$f(D) \cap M$ is homotopically incompressible in $M$.

Similarly, suppose $f(D) \cap M$ were homotopically boundary
compressible in $M$. Then there would exist an essential arc $\al$
properly embedded in $D \cap f^{-1}(M)$ such that $f(\al)$ could be
homotoped in $M$ to an arc in $\pa M$ whilst keeping its endpoints
fixed. 

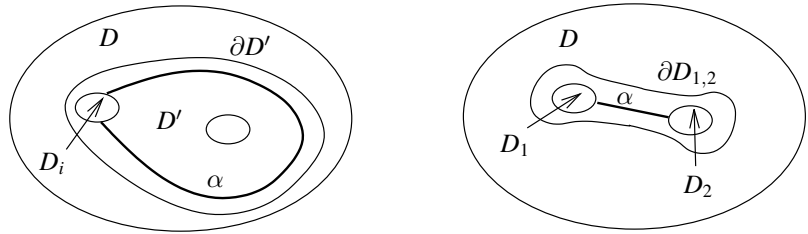
\begin{figure}[ht!]\small
\begin{center}
\begin{picture}(0,0)%
\includegraphics{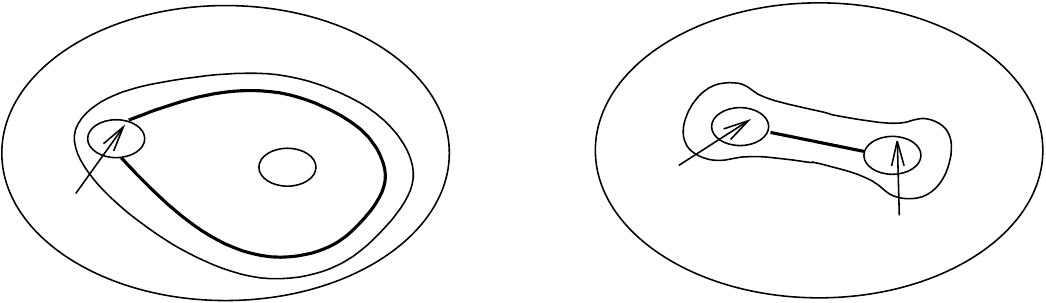}%
\end{picture}%
\setlength{\unitlength}{3947sp}%
\begin{picture}(5013,1443)(740,-1324)
\put(3818,-818){\makebox(0,0)[lb]{\smash{$D_1$}}}
\put(4968,-1077){\makebox(0,0)[lb]{\smash{$D_2$}}}
\put(1298,-151){\makebox(0,0)[lb]{\smash{$D$}}}
\put(1659,-676){\makebox(0,0)[lb]{\smash{$D'$}}}
\put(1973,-1036){\makebox(0,0)[lb]{\smash{$\alpha$}}}
\put(2116,-211){\makebox(0,0)[lb]{\smash{$\partial D'$}}}
\put(923,-931){\makebox(0,0)[lb]{\smash{$D_i$}}}
\put(4546,-518){\makebox(0,0)[lb]{\smash{$\alpha$}}}
\put(4186,-158){\makebox(0,0)[lb]{\smash{$D$}}}
\put(4801,-369){\makebox(0,0)[lb]{\smash{$\partial D_{1,2}$}}}
\end{picture}%
\caption{Two possibilities for $\al$: endpoints on the same or
  different components}
\label{fig:ensuringplanar}
\end{center}
\end{figure}

If such an arc $\al$ had both its endpoints on the same boundary
component $\pa D_i$ of $D \cap f^{-1}(M)$, it would separate $D \cap
f^{-1}(M)$ into two planar surfaces as shown in
\fullref{fig:ensuringplanar}. Here $D_i$ is a disc component of $D
\cap f^{-1}(\Hc)$. After homotoping $f(\al)$ into $\pa M$, consider a
subdisc $D' \subset D$ with image disjoint from $\Lc$ and with $f(\pa
D') \subset \pa \Hc$ running along $f(\al)$ and a subarc of $f(\pa
D_i)$, choosing this subarc so that the disc $D_i$ lies inside
$D'$. This ensures that $|D' \cap f^{-1}(\Hc)| \geq 2$.

Then $f(\pa D')$ is a curve in $\Hc$ which is homotopically trivial in
$M \cup_\phi \Hc$. By assumption, $\pi_1(\Hc)$ injects into $\pi_1(M
\cup_\phi \Hc)$, so $f(\pa D')$ is homotopically trivial in
$\Hc$. This implies that we may modify $f$ so that $f(D')$ is replaced
by a single disc in $\Hc$. But this too would reduce $|D \cap
f^{-1}(\Hc)|$, a contradiction.

Therefore we must assume that such an $\al$ has its endpoints on two
different components of $\pa (D \cap f^{-1}(M))$, say $\pa D_1$ and
$\pa D_2$ bounding discs $D_1$ and $D_2$ in $\Hc$ respectively, as
shown in the right-hand diagram of \fullref{fig:ensuringplanar}.

\begin{figure}[ht!]\small
\begin{center}
\begin{picture}(0,0)%
\includegraphics{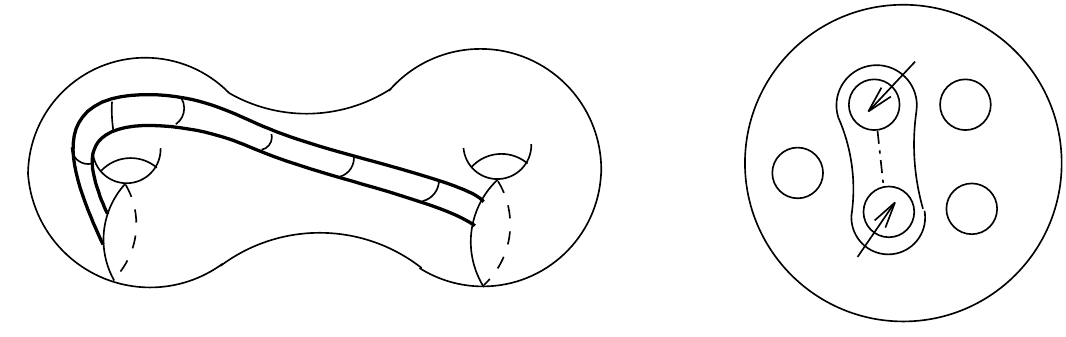}%
\end{picture}%
\setlength{\unitlength}{4144sp}%
\begin{picture}(4862,1525)(244,-996)
\put(4282,-164){\makebox(0,0)[lb]{\smash{$\alpha$}}}
\put(4429,-704){\makebox(0,0)[lb]{\smash{$D_{1,2}$}}}
\put(4315,298){\makebox(0,0)[lb]{\smash{$D_1$}}}
\put(3984,-729){\makebox(0,0)[lb]{\smash{$D_2$}}}
\put(3466,311){\makebox(0,0)[lb]{\smash{$Q$}}}
\put(2417,-943){\makebox(0,0)[lb]{\smash{$f(D_1)$}}}
\put(674,-947){\makebox(0,0)[lb]{\smash{$f(D_2)$}}}
\put(244,227){\makebox(0,0)[lb]{\smash{$H$}}}
\put(1884,-95){\makebox(0,0)[lb]{\smash{$f(\alpha)$}}}
\end{picture}%
\caption{Banding pairs of discs}
\label{fig:banding}
\end{center}
\end{figure}
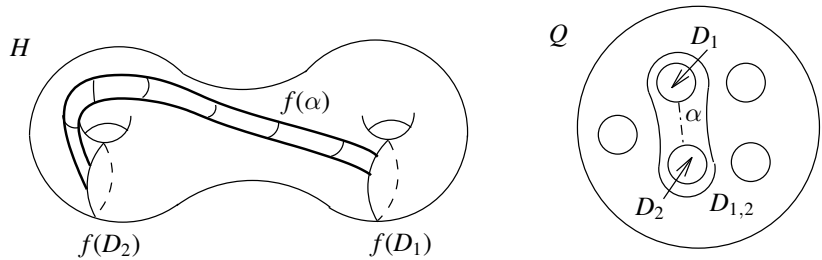

However, we may then modify $f$ by replacing both discs $D_1$ and
$D_2$ by a single disc $D_{1,2}$ consisting of $D_1$ and $D_2$ banded
together by a strip lying in a small neighbourhood of the image of the
arc $\al$ on $\pa \Hc$ pushed into $\Hc$, as in
\fullref{fig:banding}.  Note that $\pa D_{1,2}$ is embedded in $D$
although $f(\pa D_{1,2})$ is not necessarily embedded in $\Hc$. This
reduces $|D \cap f^{-1}(\Hc)|$, a final contradiction.

Thus we have shown that for any curve $\Lc$ which is homotopically
trivial in $M \cup_\phi \Hc$, we can find a spanning disc $D$ such
that the planar surface $f(D) \cap M$ is homotopically incompressible
and homotopically $\pa$--incompressible.\endproof

\proof[Proof of \fullref{thm:sfs_incomp}] By
\fullref{lem:discstructure}, we may consider a least-area immersed
planar surface $f \co Q \looparrowright M$ in the homotopy class of
$f(D) \cap M$ (\cite{Lackenby:2001}, Claim 4). Write the boundary
components of $Q$ as $q_0, \ldots, q_k$. Here $q_0$ denotes the
boundary component of $Q$ which is mapped to a geodesic curve
homotopic to $\Lc$. The other components $q_1,\ldots,q_k$ all map to
$\pa M$.

We may assume that $k \geq 2$. For, if $k=1$, the essential loop $\Lc$
is homotopic to a geodesic curve $\Lc^\ast$ in $\pa M$ bounding a disc
in $\Hc$ under the attaching map. Suppose $\widetilde{\Lc^\ast}
\subset (\pa \widetilde{M})_\ell$ is a lift of $\Lc^\ast$ to the
convex hull of $\widetilde{\Sigma}_i$, and let $A$ be an annulus
realizing the homotopy in $\widetilde{M}$ from $\widetilde{\Lc^\ast}$
to a lift $\widetilde{\Lc}$ of $\Lc$ up to $\widetilde{\Sigma}_i$.

The set $X_i \cap (\pa \widetilde{M})_\ell$ must contain some
characteristic curves, else this lift $\widetilde{\Lc}$ cannot be
homotopic to a curve in $(\pa
\widetilde{M})_\ell$. \fullref{prop:masurext} then shows that
$\widetilde{\Lc^\ast}$ intersects some characteristic curve $x \in
X_i$ lying in $(\pa \widetilde{M})_\ell$.

By an argument using Johannson's Enclosing Theorem as in
\fullref{prop:hull_equality}, some subannulus of $A$ in
$\widetilde{M}_i^\ast \subset \widetilde{M}$ can be homotoped into the
characteristic sub\-manifold $\Psi_i$ and hence into $\Psii$. But this
contradicts the assumption that $x$ is a characteristic curve. Thus $k
\geq 2$.

Consider a degree $d$ cover $R$ of $Q$ which lifts to
$\widetilde{M}$. The images of components of $\pa R$ which cover
curves $q_j$ $(j \neq 0)$ bound discs under the map $\phi \circ \pi$
and so satisfy the hypotheses of
\fullref{prop:masurext}. Recall that this ensures that such a
curve $\rho \subset (\pa \widetilde{M})_\ell$ has many quintersecting
subarcs.

Furthermore, note that $R$ has area at most $d$ times that of $Q$,
where $d$ is the degree of the cover $\pi \co \widetilde{M} \to M$.
Since $Q$ is least-area, the curvature of $Q$ inherited from $M$ is
at most $-1$ (Thm 5.5, $\S$V.A,~\cite{GHL:1987}). By Gauss--Bonnet,
\[ \mbox{Area}(R) \leq d \mbox{ Area}(Q) 
                  \leq d (2\pi(k-1)) < 2\pi kd. \] 
Recall the minimum disjoint annular width and minimum collar width
constants $\ep$ and $\de$, which depend only on $\Sigma$, $M$ and
$\widetilde{M}$. We can use these to deduce a contradiction with the
area calculation above by taking \[n \geq \frac{2\pi d}{\ep \mbox{
sinh }\de} + 2d\] as the number of quintersecting subarcs; compare
\fullref{prop:masurext}.

Consider a component $\rho$ of $\pa R$ which projects to some $q_j$
($j \neq 0$), and the $\de$--collar $N_\de(\al)$ of a quintersecting
subarc $\al$ of $\rho$. Since $\pa \widetilde{M}$ has a collar of
width $U$ in $\widetilde{M}$, and $\de \leq U$, no $\de$--collar of
another boundary component of $R$ which is a lift of some $q_j$ $(j
\neq 0)$ may intersect it.

Therefore the only boundary components of $R$ which do intersect
$N_\de(\al)$ are those which arise as lifts of $q_0$. Label such
curves $r_1,\ldots,r_t$. By construction, each $r_i$ maps into the
convex hull of some lift $\widetilde{\Sigma}_{\sg(i)}$ of $\Sigma$.

Suppose firstly that the intersection of some set $X_i$ of
distinguished curves with a boundary component $(\pa
\widetilde{M})_\ell$ contains no characteristic curves but only an
additional curve $x_\ell$. Then by our choice of $\delta$ the convex
hull $\widetilde{\Sigma}_i$ cannot intersect the $\delta$--collar of
$(\pa \widetilde{M})_\ell$. 

In particular, any curves $r_{i_j}$ lying in
$\C(\widetilde{\Sigma}_i)$ cannot intersect the $\de$--collar of
$\rho$ as it runs along $(\pa \widetilde{M})_\ell$ between
intersections with $x_\ell$. We will see below that this is enough to
deduce a contradiction. Henceforth we assume that for all $i$ and
$\ell$, $X_i \cap (\pa \widetilde{M})_\ell$ contains some
characteristic curve.

Consider an arc of intersection $y$ between a curve $r_i \subset
\C(\widetilde{\Sigma}_{\sg(i)})$ and $N_\de(\rho)$. It gives rise to a
minimal $\de$--thin arc $\xi_y$ with an endpoint $\zeta$ on
$\rho$. Suppose $\al$ is a quintersecting subarc of $\rho$ with $y
\cap N_\de(\al) \neq \emptyset$. If $y$ lies wholly in the
$\de$--collar of the interior of $\al$, $\rho$ approaches $y$ most
closely at some point in $\al$ and $\zeta \in \al$.

If instead $\zeta$ lies on a quintersecting subarc $\al'$ whose
interior is disjoint from that of $\al$, note that $y$ cannot cross
the $\de$--collar of an annulus $A_{\sg(i),j}$ parallel to a boundary
curve $x_{\sg(i),j}$ by \fullref{lem:delta}, nor the $\de$--collar
of a cut curve $x_{\sg(i),j}$ more than twice by
\fullref{lem:infinitecover} below. Assume we have decomposed $\rho$
into $n$ quintersecting subarcs; each of these intersects each
characteristic curve five times. Then $y$ can only run through the
$\de$--collar of $\al'$ and perhaps one quintersecting subarc of
$\rho$ adjacent to $\al'$.

\begin{Def}
We say that a quintersecting subarc $\al$ is \textsl{safe} if
$\bigcup_i r_i \cap N_\de(\al)$ is either connected or empty.
\end{Def}

Each minimal $\de$--thin arc corresponds to an arc $y$ which
intersects the $\de$--collar of at most two quintersecting subarcs of
$\rho$. By \fullref{lem:thinpart}, $R$ has at most $2(k-1)d$ minimal
$\de$--thin arcs. Therefore, there are at most $2(k-1)d<2kd$
quintersecting subarcs (in total counting subarcs of all $\rho \subset
\pa R$) which have two or more arcs $y$ intersecting their
$\de$--collar, and hence are not safe. 

Since there are at least $k$ boundary components in $\pa R$ which
project to some $q_j (j \neq 0)$, and each of these has $n$
quintersecting subarcs, there are at least
\[ nk-2kd=(n-2d)k \qquad \geq \qquad 
    \frac{2\pi kd}{\ep \mbox{ sinh }\de} \] safe quintersecting
subarcs. Take a safe quintersecting subarc $\al$. At most one
component $r_i \subset \pa R$ intersects $N_\de(\al)$, and it does so
in a connected subarc of $r_i$. Suppose $r_i$ maps to the convex hull
of $\widetilde{\Sigma}_{\sg(i)}$. By \fullref{prop:masurext},
$\al$ has at least five intersections with every characteristic curve
in the set $X_{\sg(i)}$.

\medskip
\textbf{Case I}\qua We deal with the simpler of two possibilities
first. Suppose that for each boundary component $(\pa
\widetilde{M})_\ell$ of $\widetilde{M}$, some curves in every $X_i$
are boundary curves of $\C(\widetilde{\Sigma}_i) \cap (\pa
\widetilde{M})_\ell$. That is, not all characteristic curves are cut
curves.

Since $\al$ is a safe quintersecting subarc, it has at least five
(hence at least three) intersections with some boundary curve
$x_{\sg(i),j} \in X_{\sg(i)}$. \fullref{lem:delta} then implies that
$\al$ runs through an annulus $A_{\sg(i),j}$ of width $\ep$ lying
outside the convex hull of $\widetilde{\Sigma}_{\sg(i)}$.

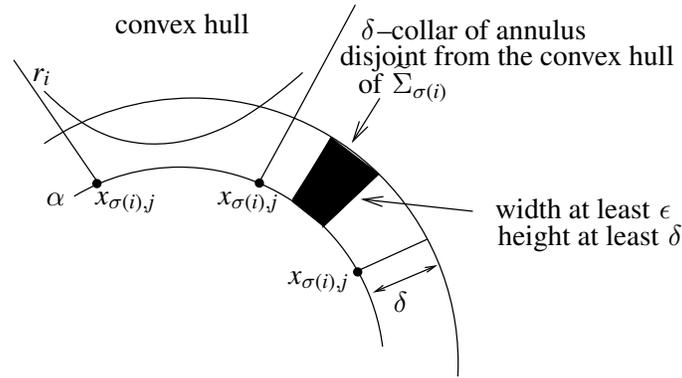
\begin{figure}[hbtp]
\begin{center}
\begin{picture}(0,0)%
\includegraphics{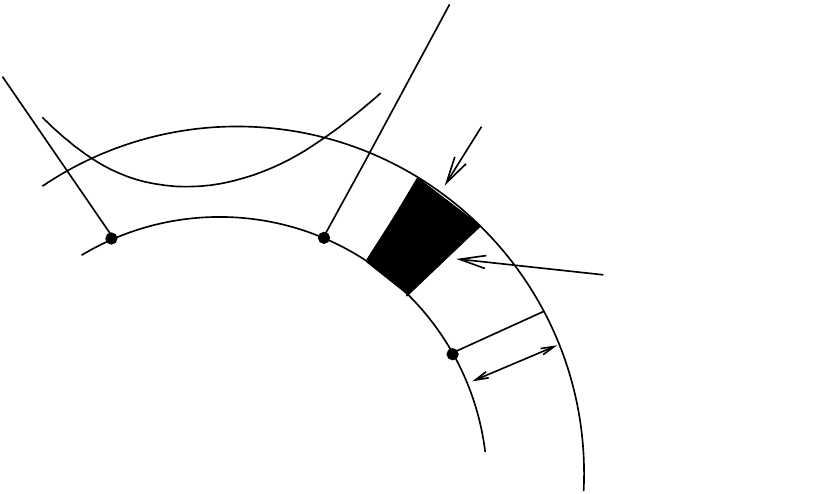}%
\end{picture}%
\setlength{\unitlength}{4144sp}%
\begin{picture}(3825,2246)(437,-1914)
\put(2401,-39){\makebox(0,0)[lb]{\smash{disjoint from the convex hull}}}
\put(3331,-961){\makebox(0,0)[lb]{\smash{width at least $\epsilon$}}}
\put(3346,-1119){\makebox(0,0)[lb]{\smash{height at least $\delta$}}}
\put(548,-136){\makebox(0,0)[lb]{\smash{$r_i$}}}
\put(2708,-1516){\makebox(0,0)[lb]{\smash{$\delta$}}}
\put(2514,112){\makebox(0,0)[lb]{\smash{$\delta$--collar of annulus}}}
\put(1661,-873){\makebox(0,0)[lb]{\smash{$x_{\sigma(i),j}$}}}
\put(927,-873){\makebox(0,0)[lb]{\smash{$x_{\sigma(i),j}$}}}
\put(2082,-1349){\makebox(0,0)[lb]{\smash{$x_{\sigma(i),j}$}}}
\put(631,-879){\makebox(0,0)[lb]{\smash{$\alpha$}}}
\put(1057,151){\makebox(0,0)[lb]{\smash{convex hull}}}
\put(2510,-211){\makebox(0,0)[lb]{\smash{of $\widetilde{\Sigma}_{\sigma(i)}$}}}
\end{picture}%
\caption{The curve $r_i$ cannot run through part of the 
         $\de$--collar of $\al \subset \rho$}
\label{fig:deltacollarcontrib}
\end{center}
\end{figure}

The image of $r_i$ in $\widetilde{M}$ cannot intersect $N_\de(\al)$
while it runs through the collar of this annulus, since the image of
$r_i$ lies in $\C(\widetilde{\Sigma}_{\sg(i)})$. By assumption, no
other boundary component of $R$ can do so either. We deduce that the
$\de$--collar of $\al$ contributes area at least $\ep \mbox{ sinh
}\de$ to $R$, as illustrated in
\fullref{fig:deltacollarcontrib}. This concludes the first case.

\medskip
\textbf{Case II}\qua Now suppose instead that $(\pa \widetilde{M})_\ell$
is a boundary component of $\widetilde{M}$ whose intersection with the
convex hull $\C(\widetilde{\Sigma}_i)$ of the embedded surface
$\widetilde{\Sigma}_i$ is all of $(\pa \widetilde{M})_\ell$. Since
this means that $\C(\widetilde{\Sigma}_i) \cap (\pa
\widetilde{M})_\ell$ is closed, it has no boundary curves. Thus for
such $i$, $X_i \cap (\pa \widetilde{M})_\ell$ consists of at least one
cut curve.

Recall \fullref{def:cutcurves}: a cut curve $x \in X_i$ is a
closed geodesic curve on some boundary component $(\pa
\widetilde{M})_\ell$ such that there exist curves $\sg^+$ and $\sg^-$
on $\widetilde{\Sigma}_i$ which are homotopic to $x$ in
$\widetilde{M}^\ast_i$ but not to each other in $\widetilde{\Sigma}_i$.

These cut curves arose when we insisted that the boundary curves of
$\Psii$ on $(\pa \widetilde{M})_\ell$ should be geodesics. Thus for
any cut curve $x$ we can find curves $\sg^+$ and $\sg^-$ in
$\widetilde{\Sigma}_i$ satisfying the definition and which were
boundary components of vertical boundary annuli for $\Psii$: see
\fullref{fig:illustrating_cut_curves}. There may be some
Seifert-fibred solid torus components of the characteristic
submanifold $\Psi_i$ lying between the components of $\Psii$, but they
will make no difference to our argument.

\begin{Lem} 
\label{lem:infinitecover}
Suppose $r_i \subset \pa R$ maps to the convex hull
$\C(\widetilde{\Sigma}_{\sg(i)})$, and suppose that $y$ is a connected
subarc of $r_i$ lying in the $\de$--collar of $\al$. Then this arc $y$
intersects the $\de$--collar of a cut curve $x_{\sg(i),j}$ at most
twice.
\end{Lem}

\proof Suppose not, so $y$ intersects the $\de$--collar of the
cut curve $x_{\sg(i),j}$ three times. Choose one of the endpoints of $y$
as a basepoint for $\widetilde{M}$. Consider the cover
$\widetilde{M}_\infty \to \widetilde{M}$ such that
$\pi_1(\widetilde{M}_\infty)= \pi_1(\widetilde{\Sigma}_{\sg(i)})$, and
give $\widetilde{M}_\infty$ a basepoint which projects to the
basepoint for $\widetilde{M}$. Note that the manifold
$\widetilde{M}_\infty$ is homeomorphic to $\widetilde{\Sigma}_{\sg(i)}
\times I$, and moreover that $\widetilde{\Sigma}_{\sg(i)}$ lifts
homeomorphically to $\widetilde{M}_\infty$. Write
$\widetilde{\Sigma}_{\sg(i),\infty}$ for such a lift.

Take a $\de$--collar of the boundary of $\widetilde{M}_\infty$, and
denote it by $N_\de(\pa \widetilde{M}_\infty)$. By decreasing $\de>0$
again if necessary (see before \fullref{def:collarwidth}), the
intersection of the convex hull of
$\widetilde{\Sigma}_{\sg(i),\infty}$ with this $\de$--collar is a
collar on the intersection of the convex hull with the boundary of
$\widetilde{M}_\infty$. In this way we obtain
\[ \C(\widetilde{\Sigma}_{\sg(i),\infty}) \cap N_\de(\pa
\widetilde{M}_\infty) = \left
(\C(\widetilde{\Sigma}_{\sg(i),\infty}) \cap \pa \widetilde{M}_\infty
\right ) \times I. \]
Since the cover $\widetilde{M}_\infty$ only contains loops which lie
in $\pi_1(\widetilde{\Sigma}_{\sg(i)})$, the convex hull has been cut
along every cut curve in $X_{\sg(i)}$. Thus every component of the
above intersection is the lift of a collar on the appropriate
intersection of the original characteristic submanifold
$\Psi_{\sg(i)}$ with $\pa \widetilde{M}$, where we assume that
boundary curves on $\pa \widetilde{M}$ remain disjoint rather than
being made geodesic.

The loop $r_i$ is homotopic to a loop in
$\widetilde{\Sigma}_{\sg(i)}$, and passes through the basepoint of
$\widetilde{M}$. When we lift $r_i$ to $\widetilde{M}_\infty$, it
remains a loop based at the chosen lift of the basepoint, lying in the
convex hull $\C(\widetilde{\Sigma}_{\sg(i),\infty})$. The subarc $y$
lifts to a subarc in $\widetilde{M}_\infty$ lying entirely in one
component of the intersection of
$\C(\widetilde{\Sigma}_{\sg(i),\infty})$ with $N_\de(\pa
\widetilde{M}_\infty)$. Note also that the cut curve $x_{\sg(i),j}$ is
homotopic to a loop in $\widetilde{\Sigma}_{\sg(i)}$ so that it can be
lifted to a loop in $\widetilde{M}_\infty$ as well.

\begin{figure}[ht!]\small
\begin{center}
\begin{picture}(0,0)%
\includegraphics{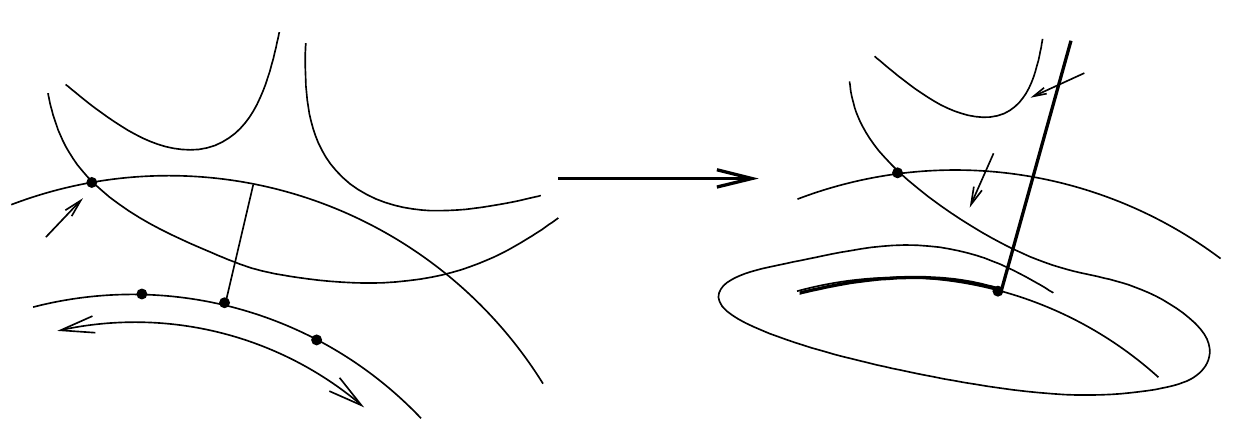}%
\end{picture}%
\setlength{\unitlength}{4144sp}%
\begin{picture}(5740,1920)(16,-1159)
\put(5003,417){\makebox(0,0)[lb]{\smash{convex hull}}}
\put( 68,-451){\makebox(0,0)[lb]{\smash{basepoint}}}
\put( 16,-810){\makebox(0,0)[lb]{\smash{$\rho$}}}
\put(2341,-564){\makebox(0,0)[lb]{\smash{$r_i$}}}
\put(1087,-600){\makebox(0,0)[lb]{\smash{$x_{\sigma(i),j}$}}}
\put(1059,-968){\makebox(0,0)[lb]{\smash{$\alpha$}}}
\put(736,-225){\makebox(0,0)[lb]{\smash{$y$}}}
\put(3923,-751){\makebox(0,0)[lb]{\smash{lift of $x_{\sigma(i),j}$}}}
\put(2581, 73){\makebox(0,0)[lb]{\smash{lift to $\widetilde{M}_\infty$}}}
\put(4170,105){\makebox(0,0)[lb]{\smash{lift of $y$}}}
\put(3578,629){\makebox(0,0)[lb]{\smash{$\widetilde{\Sigma}_{\sigma(i),
\infty}$}}}
\put(278,607){\makebox(0,0)[lb]{\smash{$\widetilde{\Sigma}_{\sigma(i)}$}}}
\end{picture}%
\caption{The subarc $y$ of $r_i$ and its lift to
$\widetilde{M}_\infty$}
\label{fig:infinitecover}
\end{center}
\end{figure}
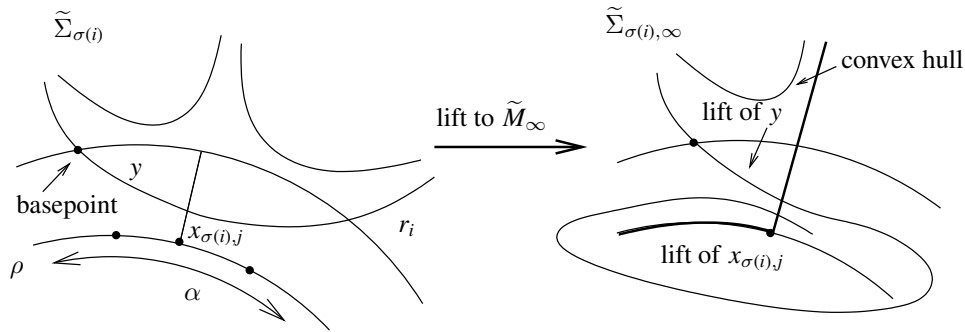

Since $y$ crosses the $\de$--collar of a cut curve $x_{\sg(i),j}$
three times, we can choose two of these intersections which have
coherent orientations. Thus there is a based loop (see
\fullref{fig:infinitecover}) lying entirely in one component of
$\C(\widetilde{\Sigma}_{\sg(i),\infty}) \cap N_\de(\pa
\widetilde{M}_\infty)$ which has non-zero geometric intersection with
the $\de$--collar of some lift $\widetilde{x}_{\sg(i),j}$ of the cut
curve. But this component is a lift of a component of
\[ \left ( \Psi_{\sg(i)} \cap \pa \widetilde{M} \right ) \times I \]
and so cannot contain such a based loop crossing the lift of any cut
curve. This is a contradiction, proving
\fullref{lem:infinitecover}.\endproof

Therefore a curve $r_i$ can intersect the $\de$--collar of at most two
intersections of $\al$ with $x_{\sg(i),j}$. Since
\fullref{prop:masurext} ensures that $\al$ intersects
$x_{\sg(i),j}$ at least five times, there are two consecutive
intersections between which $r_i$ does not intersect the $\de$--collar
of $\al$. By choice of $\ep$, these intersections are at least
distance $\ep$ apart. This ensures that $\al$ contributes at least
$\ep \mbox{ sinh }\de$ to the area of $R$ between these two
intersections, concluding the second case.

We have shown the existence of at least $2\pi kd/(\ep \mbox{ sinh
}\de)$ subarcs of $\pa R$ each contributing area at least $\ep \mbox{
sinh }\de$ to $R$. This implies that
\[ \mbox{Area}(R) \geq \frac{2\pi kd}{\ep \mbox{ sinh }\de}
\times \ep \mbox{ sinh } \de \geq 2 \pi kd \] which contradicts the
Gauss--Bonnet calculation above. Thus no such surface $R$ can exist,
and our original assumptions cannot hold. In particular, no essential
loop in $\Sigma$ bounds a disc $f \co D \looparrowright M \cup_\phi \Hc$. 

Let us summarize. For a suitable $\pi_1$--injective virtually embedded
surface $\Sigma$, there exists a constant $n$ as described above
depending only on $\Sigma$, $M$ and $\widetilde{M}$. Choosing suitable
pseudo-Anosov maps $h_{\ell} \co \pa H_{\ell} \to \pa H_{\ell}$ and
homeomorphisms $\phi'_{\ell} \co (\pa M)_\ell \to \pa H_\ell$, there
exist integers $(N_\ell)_{\textrm{min}}$ depending on $n$ and which
have the following property. Whenever $N_{\ell} \geq
(N_\ell)_{\textrm{min}}$ for all $\ell=1,\ldots,m$, the surface
$\Sigma$ remains $\pi_1$--injective in $M \cup_\phi \Hc$ under the
handlebody attachments $\phi_{\ell}=h_{\ell}^{N_{\ell}} \circ
\phi_{\ell}' \co (\pa M)_{\ell} \to \pa H_{\ell}$. This proves
\fullref{thm:sfs_incomp}.\endproof

\proof[Proof of \fullref{thm:sfsubgp}] By work of
Cooper--Long--Reid~\cite{CooperLongReid:1997}, any 3--manifold $M$
satisfying the hypotheses of \fullref{thm:sfsubgp} contains a
surface $\Sigma$ satisfying the assumptions of
\fullref{thm:sfs_incomp}. Indeed, there may be infinitely many
such surfaces. They are formed by connecting two copies of an
essential surface realizing some boundary slope on $\pa M$ by
high-genus tubes wrapping many times round the appropriate components
of $\pa M$, and then compressing if necessary.

\fullref{thm:sfs_incomp} shows that such surfaces $\Sigma$ remain
$\pi_1$--injective under appropriate gluing maps $\phi$. In
particular, $\pi_1(M \cup_\phi \Hc)$ contains $\pi_1(\Sigma)$ as a
surface subgroup. This enables us to deduce
\fullref{thm:sfsubgp}.\endproof

\section{Proof of main theorem: general case}
\label{sec:many_components}

In \fullref{sec:incomp} we proved our main theorem in the case
where $M$ does not have cusps. We now show that our arguments
also extend to the cusped case.

\begin{Thm}
\label{thm:sfsubgp_many}
Suppose $M \neq B^3$ is a simple 3--manifold with $m+m' \geq 1$
boundary components, $m'$ of which are tori. Let
$\Hc=\{H_1,\ldots,H_{m+m'}\}$ be a collection of handlebodies and
solid tori whose genera match those of $\pa M$. Let $M \cup_\phi \Hc$
denote the closed 3--manifold obtained by gluing each boundary
component $(\pa M)_{\ell}$ to $\pa H_{\ell}$ by a homeomorphism
$\phi_\ell$.

Suppose moreover that $h_\ell \co \pa H_\ell \ra \pa H_\ell$ is a
homeomorphism which is either a pseudo-Anosov homeomorphism whose
stable lamination is of full type, or an Anosov homeomorphism,
according to whether $H_\ell$ has genus at least two or is a solid
torus respectively.

Given homeomorphisms $\phi'_\ell \co (\pa M)_\ell \to \pa H_\ell$,
there exist integers $(N_\ell)_{\textrm{min}}$ such that if $\phi_\ell
= h_\ell^{N_\ell} \circ \phi'_\ell$ with $N_\ell \geq
(N_\ell)_{\textrm{min}}$ for all $\ell$, the group $\pi_1(M
\cup_{\phi} \Hc)$ contains a surface subgroup.
\end{Thm}

Recall from~\cite{CassonBleiler:1988} that an Anosov homeomorphism of
a torus is a map whose trace under any geometric
$SL(2,\Zb)$--representation has absolute value strictly greater than
two. They fulfil a similar dynamical purpose to pseudo-Anosov maps.

\fullref{thm:sfsubgp_many} may be deduced from the appropriate
analogue of \fullref{thm:sfs_incomp} as in the non-cusped
case. Since $M$ is Haken, we may give $M$ (minus its toral boundary) a
complete hyperbolic structure with totally geodesic boundary.

\begin{Thm}
\label{thm:sfs_incomp_many}
Let $M$ be a hyperbolic 3--manifold with $m$ totally geodesic boundary
components and $m'$ cusps, $m+m' \geq 1$. Let $\Sigma$ be a connected,
orientable, closed, immersed surface in $M$ with a finite cover which
lifts to a non-peripheral embedded incompressible surface in a finite
regular cover $\pi \co \widetilde{M} \to M$ of degree $d$.

As in \fullref{thm:sfsubgp_many}, let $M \cup_\phi \Hc$ denote the
closed 3--manifold obtained by filling each of the boundary components
or cusps of $M$ via homeomorphisms $\phi_\ell \co (\pa M)_{\ell} \ra
\pa H_{\ell}$, and consider maps $h_\ell$ satisfying the condition
given in that theorem.

Given homeomorphisms $\phi'_\ell \co (\pa M)_\ell \to \pa H_\ell$,
there exist integers $(N_\ell)_{\textrm{min}}$ such that if $\phi_\ell
= h_\ell^{N_\ell} \circ \phi'_\ell$ with $N_\ell \geq
(N_\ell)_{\textrm{min}}$ for all $\ell$, the surface $\Sigma$ remains
$\pi_1$--injective in the resulting manifold $M \cup_\phi \Hc$.
\end{Thm} 

\proof[Proof of \fullref{thm:sfs_incomp_many}] We argue by contradiction
as in \fullref{thm:sfs_incomp}. Apply
\fullref{lem:discstructure} to a disc $f \co D \looparrowright M
\cup_\phi \Hc$ spanning a loop $\Lc$ in $\Sigma$, and obtain once more
a least-area planar surface $f \co Q \looparrowright M$ with $k+1$
boundary components and $k'$ punctures. We say that the boundary
components mapping to totally geodesic boundary components of
$\widetilde{M}$ are \textsl{non-toral curves}, and label them
$q_1,\ldots,q_k$. We also label the \textsl{punctures}
$p_1,\ldots,p_{k'}$. Take a cover $\tilde{f} \co R \looparrowright
\widetilde{M}$ as before.

By the same arguments as for the non-cusped case, we can ensure $k+k'
\geq 2$. There are then two possible cases: either $\Sigma$ contains
accidental parabolics or it does not. Since $\Sigma$ is closed, it
contains no genuine parabolics.

Before we examine these cases separately, we need to generalize the
proof of one of our results to the cusped case. Recall the proof of
\fullref{prop:hull_equality}: at one point we used the
assumption that $\pi_1(M)$ had no parabolic elements. With this
assumption we deduced that $\Psii \cap \pa \widetilde{M} =
\C(\widetilde{\Sigma}_i) \cap \pa \widetilde{M}$.

In fact the same proof works in the cusped case, but we need to
justify this. Note that $M$ still cannot be virtually fibred over
$\Sigma$ since $\Sigma$ is closed. The only potential problem is the
set $P$ of parabolic fixed points which occurs in the identity
(Theorem 3.14 of~\cite{MatsuTani:1998})
\[ \Lambda(\Gamma_1) \cap \Lambda(\Gamma_2) =
   \Lambda(\Gamma_1 \cap \Gamma_2) \cup P \] for two geometrically
finite subgroups $\Gm_1,\Gm_2$ of a Kleinian group $\Gm$. 

The set $P$ is defined in~\cite{MatsuTani:1998} as the set of points
$\zeta$ in the complement $\Omega(\Gm_1 \cap \Gm_2)$ of $\Lambda(\Gm_1
\cap \Gm_2)$ such that stab$_{\Gm_1}(\zeta)$ and stab$_{\Gm_2}(\zeta)$
generate a parabolic abelian group of rank 2 and stab$_{\Gm_1}(\zeta)$
$\cap$ stab$_{\Gm_2}(\zeta) = \{\mbox{id}\}$.

The set $P$ is only non-empty when the appropriate stabilizers of
$\zeta \in P$ are both of rank 1 and intersect trivially. In our case,
$\Gm_1$ is the fundamental group of some lift $\widetilde{\Sigma}_i$.
The group $\Gm_2$ is the fundamental group of some boundary component
$(\pa \widetilde{M})_\ell$ which is either totally geodesic or
cuspidal.

If it is cuspidal, the group $\Gm_2$ contains a rank 2 stabilizer of
$\zeta$ if and only if $(\pa \widetilde{M})_\ell$ is the appropriate
toral boundary component for the cusp corresponding to $\zeta$. If
$(\pa \widetilde{M})_\ell$ is non-toral or if it corresponds to a
different cusp, only the identity element stabilizes the parabolic
fixed point $\zeta$. The group $\mbox{stab}_{\Gm_2}(\zeta)$ can never
be of rank 1. In our application, therefore, $P$ is always empty.

We deduce that \fullref{prop:hull_equality} still holds when
$M$ has cusps. Thus, assuming that the characteristic submanifold
$\Psii$ has geodesic boundary on $\pa \widetilde{M}$, the convex hull
of $\widetilde{\Sigma}_i$ intersects each torus $(\pa
\widetilde{M})_\ell$ in immersed annuli whose boundary curves are
characteristic curves of $\Psii$ and whose image may be all of $(\pa
\widetilde{M})_\ell$.

Now we can return to the main argument. By Gauss--Bonnet applied to
$Q$, 
\[ \mbox{Area}(R) \leq d \mbox{ Area}(Q) \leq d (2\pi(k+k'-1)) 
< 2\pi (k+k')d. \] Choose disjoint horoball neighbourhoods of each
cusp in $\widetilde{M}$ which are disjoint from the
$U$\negthinspace--collars of every non-toral curve. Suppose moreover
that the $\de$--collar of each horoball neighbourhood is disjoint from
any convex hull which is disjoint from the toral boundary
corresponding to this cusp.

By homotopic $\pa$--incompressibility, each puncture $\tilde{p}$ in
$R$ covering a puncture in $Q$ has a half-open annulus neighbourhood
$A(\tilde{p})$ whose image lies in one of these horoball
neighbourhoods. The boundary of each horoball neighbourhood is a torus
$T_\ell$ which inherits a Euclidean path metric from the metric on
$\widetilde{M}$.

\medskip
\textbf{Case A: no accidental parabolics}\qua Suppose $\Sigma$ contains
no accidental parabolics: non-trivial loops in $\Sigma$ homotopic to
loops in the cusp boundary of $M$. Then, for any lift
$\widetilde{\Sigma}_i$, the convex hull $\C(\widetilde{\Sigma}_i)$
avoids the cusps of $\widetilde{M}$.

Write $\rho(\tilde{p})$ for the boundary component of $A(\tilde{p})$
on some $T=T_\ell$. We say that such curves are \textsl{toral curves}
in $R$. By the appropriate and much easier analogue of
\fullref{prop:masurext} for these toral attachments, we may
ensure that each toral curve $\rho(\tilde{p})$ has uniformly large
geometric intersection number with some additional curve $x_\ell$. As
these intersections occur at least distance $\ep$ apart, we may also
ensure that all curves $\rho(\tilde{p})$ are longer than any given
fixed length.

The area of $A(\tilde{p})$ is at least the length of $\rho(\tilde{p})$
(see e.g.~\cite{HassRubWang:1999}, proof of Thm 4.3). By the above
argument we can make all toral curves at least $2d(\pi + \ep \mbox{
sinh }\de)$ long so that each annulus $A(\tilde{p})$ contributes
$2d(\pi + \ep \mbox{ sinh }\de)$ to the area of $R$. There are at
least $k'$ such annuli.

As in the proof of \fullref{thm:sfs_incomp}, we can subdivide
every non-toral curve into at least $n \geq 2\pi d/(\ep \mbox{ sinh
}\de) + 2d$ quintersecting subarcs. Over all non-toral curves, all but
at most $2(k+k')d$ of these subarcs contribute area $\ep \mbox{ sinh
}\de$ to $R$. Thus
\begin{eqnarray*}
 \mbox{Area}(R) & \geq & \left ( \frac{2\pi k d}{\ep \mbox{ sinh }\de} +
2kd - 2(k+k')d \right ) \times \ep \mbox{ sinh }\de 
    + 2k'd (\pi + \ep \mbox{ sinh } \de )  
\end{eqnarray*}
and so $\mbox{Area}(R) \geq 2\pi(k+k')d$, which contradicts the
Gauss--Bonnet calculation.

In fact, we could have obtained a contradiction by showing that each
pair of intersections of $\rho(\tilde{p})$ with $x_\ell$ picked up a
quantum of area \emph{outside} $A(\tilde{p})$. We need to use this
idea in the general case.

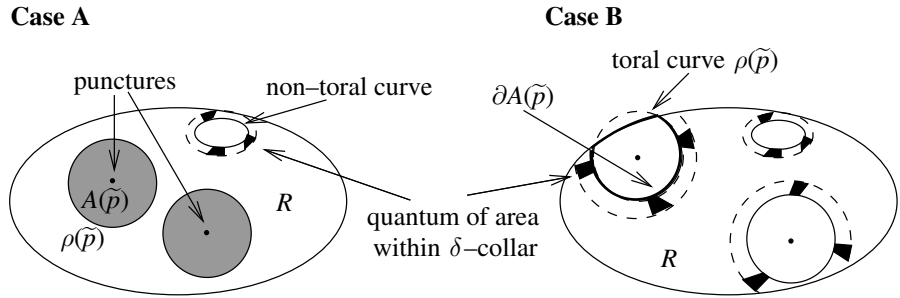
\begin{figure}[ht!]\small
\begin{center}
\begin{picture}(0,0)%
\includegraphics{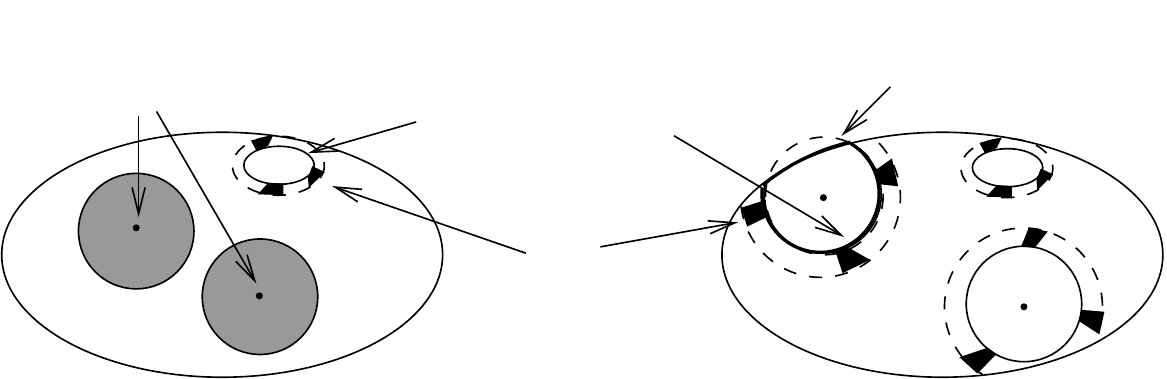}%
\end{picture}%
\setlength{\unitlength}{4144sp}%
\begin{picture}(5325,1730)(123,-933)
\put(3003,214){\makebox(0,0)[lb]{\smash{$\partial A(\tilde{p})$}}}
\put(3331,689){\makebox(0,0)[lb]{\smash{\textbf{Case B}}}}
\put(1665,270){\makebox(0,0)[lb]{\smash{non--toral curve}}}
\put(2270,-510){\makebox(0,0)[lb]{\smash{quantum of area}}}
\put(514,306){\makebox(0,0)[lb]{\smash{punctures}}}
\put(552,-408){\makebox(0,0)[lb]{\smash{$A(\tilde{p})$}}}
\put(2305,-695){\makebox(0,0)[lb]{\smash{within $\delta$--collar}}}
\put(422,-637){\makebox(0,0)[lb]{\smash{$\rho(\tilde{p})$}}}
\put(1714,-428){\makebox(0,0)[lb]{\smash{$R$}}}
\put(4015,-758){\makebox(0,0)[lb]{\smash{$R$}}}
\put(3733,432){\makebox(0,0)[lb]{\smash{toral curve $\rho(\tilde{p})$}}}
\put(136,689){\makebox(0,0)[lb]{\smash{\textbf{Case A}}}}
\end{picture}%
\caption{Comparing the area calculations in Case A and Case B}
\label{fig:cuspsCaseA}
\end{center}
\end{figure}

\medskip
\textbf{Case B: accidental parabolics}\qua Now suppose that the surface
$\Sigma$ contains some accidental parabolics. Consider a boundary
component $T$ of a horoball neighbourhood of a cusp in $\widetilde{M}$
which has non-empty intersection with some $\C(\widetilde{\Sigma}_i)$,
and suppose that $A(\tilde{p})$ lies in the corresponding horoball
neighbourhood.

In this case, the boundary curve $\pa A(\tilde{p})$ may not lie
entirely in $T$. Define a \textsl{toral curve} $\rho(\tilde{p})$
instead to be a curve freely homotopic to a core curve of
$A(\tilde{p})$ lying on $T$. Also define $N_\de(A(\tilde{p}))$ to be
the set of all points in $R$ which lie within $\de$ of
$A(\tilde{p})$. There are two cases.

\medskip
\textbf{Case B(I)}\qua $\C(\widetilde{\Sigma}_i) \cap T \neq T$ for all
$i$\qua By construction, the only boundary components of $R$ which
intersect $N_\de(A(\tilde{p}))$ are those components $r_i$ which cover
$q_0$ and which map to some $\C(\widetilde{\Sigma}_{\sg(i)})$
intersecting $T$.

Since $\C(\widetilde{\Sigma}_{\sg(i)}) \cap T$ is never all of $T$, it
has some boundary curve $x'_{\sg(i),j}$ on $T$. Again we may adapt
\fullref{prop:masurext} to the case of these toral
attachments, since they are the composition of a large power of an
Anosov map with the map $\phi'$. Thus $\rho(\tilde{p})$ can be
subdivided into many quintersecting subarcs, each intersecting
$x'_{\sg(i),j}$ at least five times. But none of the intersections can
occur while $r_i$ runs through $A(\tilde{p})$, so we can replace
$\rho(\tilde{p})$ by $\pa A(\tilde{p})$ in this statement.

Taking the powers of the Anosov map in the solid torus attachments to
be sufficiently large, we assume that each curve $\pa
A(\tilde{p})$ has been subdivided into at least $n \geq 2\pi d/(\ep
\mbox{ sinh }\de) + 2d$ quintersecting subarcs.

By analogy with \fullref{def:thinpart}, define minimal
$\de$--thin arcs in $R$ to be possibly degenerate minimal arcs of
length at most $\de$ with endpoints on $\pa A(\tilde{p})$ and some
$r_i$ respectively. Such an arc may be collapsed to a single
point. Adapting \fullref{lem:thinpart} to this case, $R$ has $\leq
2(k+k')d$ minimal $\de$--thin arcs. Suppose that $\kappa$ of these
occur within the $\de$--collar of a non-toral curve.

This leaves at most $2(k+k')d-\kappa$ subarcs of curves $\pa
A(\tilde{p})$ which have two or more curves $r_i$ intersecting their
$\de$--collar. By the appropriate version of \fullref{lem:delta},
every other subarc is safe and contributes area $\ep \mbox{ sinh }\de$
to $R$.

Subdivide each non-toral curve into at least $n \geq 2\pi d/(\ep
\mbox{ sinh }\de) +2d$ subarcs. Over all non-toral curves, at most
$\kappa$ subarcs fail to contribute area $\ep \mbox{ sinh
}\de$. Combining these contributions,
\[ \frac{\mbox{Area}(R)}{\ep \mbox{ sinh }\de} 
   \geq \left ( \frac{2 \pi kd}{\ep \mbox{ sinh }\de} + 2kd - \kappa
\right ) + \left ( \frac{2 \pi k'd}{\ep \mbox{ sinh }\de} + 2k'd -
(2(k+k')d - \kappa) \right ) \]
and hence we obtain $\mbox{Area}(R) \geq 2\pi(k+k')d$, contradicting
Gauss--Bonnet. 

\medskip
\textbf{Case B(II)}\qua $\C(\widetilde{\Sigma}_i) \cap T = T$ for some
$i$\qua We may combine the arguments above with those explained in Case
II of the proof of \fullref{thm:sfs_incomp}.

This proves \fullref{thm:sfs_incomp_many} and hence
\fullref{thm:sfsubgp_many}.\endproof

\bibliographystyle{gtart}
\bibliography{link}

\end{document}